\documentclass{article}[12pt]
\usepackage{amsmath,amsfonts,amsthm}
\usepackage{amssymb}

\usepackage{epsfig}

\usepackage{graphics}
\usepackage{graphicx}

\usepackage{latexsym}
\usepackage{graphics}
\usepackage{amsmath}
\usepackage{amsthm}
\usepackage{xspace}
\usepackage{amssymb}
\usepackage{color}
\usepackage{cite}
\usepackage{latexsym}
\usepackage{graphics}
\usepackage{amsmath}
\usepackage{amsthm}
\usepackage{xspace}
\usepackage{amssymb}
\usepackage{epsfig}

\usepackage{enumitem}

\newcommand{\beql}[1]{\begin{equation}\label{#1}}
\newcommand{\eeql}{\end{equation}}
\newcommand{\eqn}[1]{(\ref{#1})}

\newcommand{\sasha}[1]{#1}

\newcommand{\R}{\mathbb{R}}
\newcommand{\pr}{\mathbb{P}}
\newcommand{\E}{\mathbb{E}}

\def\P{{\mathbb P}}

\newcommand{\cx}{{\cal X}}

\newcommand{\cj}{{\cal J}}

\newcommand{\bI}{{\bf I}}

\newtheorem{thm}{Theorem}[section]
\newtheorem{lem}[thm]{Lemma}
\newtheorem{prop}[thm]{Proposition}
\newtheorem{cor}[thm]{Corollary}
\newtheorem{assumption}[thm]{Assumption}
\newtheorem{definition}[thm]{Definition}
\newtheorem{condition}[thm]{Condition}
\newtheorem{remark}[thm]{Remark}
\newtheorem{example}[thm]{Example}
\begin{document}

\title{Large-scale distributed synchronization systems, using a 
cancel-on-completion redundancy mechanism  
%\\(DRAFT)
}

\author
{
Alexander L. Stolyar \\
The Grainger College of Engineering\\
ISE Department
and Coordinated Science Lab\\
University of Illinois at Urbana-Champaign\\
Urbana, IL 61801, USA \\
\texttt{stolyar@illinois.edu}
}

\date{\today}

\maketitle

\begin{abstract}

We consider a class of multi-agent distributed synchronization systems, which are modeled as $n$ particles moving on the real line.
This class generalizes the model of a multi-server queueing system, considered in \cite{St2021-coc},
employing so-called cancel-on-completion (c.o.c.) redundancy mechanism, but is motivated by other applications as well.
In the multi-server queueing system a particle location represents a server workload. Under c.o.c. mechanism, 
when a job of class $j$ arrives, it selects $d_j$ particles uniformly at random, which try to jump forward, by random distances, but their advance is truncated at the new location of the $k_j$-th left-most selected particle ($k_j \le d_j$). Between jumps all particles move to the left at constant speed, but cannot cross point $0$ (workload cannot be less than $0$). Thus, the multi-server queueing system is modeled as a particle system, regulated at the left
boundary point. 

 The more general model of this paper is such that particles evolve the same way as the in left-regulated system, 
 but we allow regulation boundaries on either side, or both sides, or no regulation at all.
 We consider the mean-field asymptotic regime, when the number of particles $n$ and the job arrival rates go to infinity, while the job arrival rates 
per particle remain constant. The system state for a given $n$ is the empirical distribution of the particles' locations. Our results include: the existence/uniqueness of fixed points of mean-field limits (ML), which describe the limiting dynamics of the system; conditions for the steady-state asymptotic independence (concentration, as $n \to\infty$, of the stationary distribution on a single state, which is necessarily an ML fixed point); the limits, as $n \to\infty$, of the average velocity at which unregulated (free) particle system advances.
In particular, our results for the left-regulated system unify and generalize the corresponding results in  \cite{St2021-coc}.
Our technical development is such that the systems with different types of regulation are analyzed within a unified framework.
In particular, these systems are used as tools for analysis of each other.

\end{abstract}


\iffalse
\noindent
{\em Abbreviated Title:} 
\fi


%\newpage

\section{Introduction}
\label{sec-intro}

\subsection{Model introduction and motivation}

We consider a class of particle systems which generalizes the class of multi-server queueing models, considered in \cite{St2021-coc},
with so-called cancel-on-completion (c.o.c.) redundancy. The evolution of server workloads in a multi-server system with c.o.c. redundancy 
can be viewed as evolution of particles on the real line, with regulation at the left boundary at $0$ -- particles (server workloads) are not allowed 
to move to the left of $0$. The more general model of this paper is such that, in addition to left-regulated particle systems,
we also consider right- and/or left-regulated, as well as unregulated (free) systems.
We study the mean-field asymptotic regime, when the number of particles and the job arrival rates go to infinity, while the job arrival rates 
per particle remain constant. We obtain results on both the asymptotics of stationary distributions and on the behavior of mean-field limits, which describe the limiting dynamics of the system. (Unregulated systems were considered already in \cite{St2021-coc}, but not their asymptotic behavior.) 
A feature of our technical development is that the systems with different types of regulation are analyzed within a unified framework.
In particular, such systems are used as tools for analysis of each other.

A more detailed description of the multi-server queueing model with c.o.c. redundancy is as follows. There are $n$ identical servers, processing work at unit rate. New jobs arrive as a Poisson process of rate $\lambda^n n$. There is a finite set of job classes, and an arriving job is of class $j$ with probability $\pi_j^n$.
A class $j$ job consists of $d_j$ components, which are placed on $d_j$ servers selected uniformly at random;
the components' sizes (service times) $(\xi^{(j)}_1,\ldots, \xi^{(j)}_{d_j})$ are drawn according to an exchangeable distribution $F_j$.
Each server processes its work (components of different jobs) in the First-Come-First-Served order. 
A class $j$ job service is complete as soon as any $k_j$ components (with $k_j \le d_j$) of the job complete their service, at which point the unfinished service of the remaining $d_j -k_j$ components of the job is ``canceled,'' and they are immediately removed from the system;
we call this $(d_j,k_j)$-c.o.c. redundancy.
A server workload at a given time $t$ is its current unfinished work, namely the amount of time 
it will take the server to finish all work, assuming no new jobs will arrive in the system after $t$.

One can view the process described above as a particle system, if each server $i$ is associated with a ``particle,'' and its workload $W_i$
is interpreted as the particle ``location'' on the real line. With this view, 
an arriving class-$j$ ``job'' selects $d_j$ particles uniformly at random,
and their random non-negative potential jump sizes $(\xi^{(j)}_1,\ldots, \xi^{(j)}_{d_j})$ are drawn from distribution $F_j$. So, if a selected particle $i$ current location is $W_i$, its potential new location is $W_i + \xi_i^{(j)}$. However, the actual new location particle $i$ jumps to is ``truncated'' at point $W^*$, which is the $k_j$-th smallest potential new location among the selected particles (due to the c.o.c. mechanism), so that the particle actually jumps to 
$\max\{W_i, \min\{W_i + \xi_i^{(j)}, W^*\}\}$.
Between ``job'' arrivals, each particle moves to the ``left'' at constant speed $v=1$,
but cannot move to the left of the ``regulation'' boundary at $0$; a particle system with this feature will be referred to as {\em left-regulated.} 
Due to the jump truncation feature (c.o.c. mechanism) the system is ``non-work-conserving'' in the sense that the expected total displacement of the selected particles after a class-$j$ job arrival
is not constant -- it, in general, depends on the system state at the time of arrival, namely on the relative locations of selected particles.
The non-work-conservation is the model feature which substantially complicates the system analysis (see section 1.1 in \cite{St2021-coc}).

We study the mean-field asymptotic regime, when $n\to\infty$, $\lambda^n \to \lambda$, and $\pi^n_j \to \pi_j$ for all $j$. 
The system state for a given $n$ is the empirical distribution of particles' locations.
In this regime, deterministic trajectories arising as limits of the process dynamics are called mean-field limits (ML).
We say that the steady-state asymptotic independence (SSAI) holds when 
the sequence of stationary distributions converges to the distribution concentrated on a single state (which is then necessarily an ML fixed point).
Some of the main 
questions 
in this asymptotic regime concern the existence and characterization of ML fixed points and conditions under which the SSAI holds.
Paper \cite{St2021-coc} focused on these questions for  the left-regulated system.

In this paper we study particle systems with the same job arrival and jump truncation (c.o.c.) mechanism, under the same asymptotic regime, 
but more general in that we allow systems regulated on the left and/or right, as well as unregulated (free) systems. 
It is also convenient to have the speed $v \ge 0$, at which each particle moves to the left, as an additional parameter,
and consider the dependence of the system properties on $v$, for fixed parameters $\lambda$, $\{\pi_j\}$, and $\{F_j\}$.

Additional question of interest for a free system is: If $v_n$ is the average 
velocity at which the entire particle system (say, its mean) advances to the ``right,'' does $\lim_{n\to\infty} v_n$ exist? If so, is the limit determined by the speed parameter of the free system ML fixed points (if any).
Note that in a free system there is the obvious relation between ML fixed points and MLs that are traveling waves. Namely, for any $u \ge 0$,
an ML fixed point for speed parameter $v$ corresponds to a traveling wave ML moving right with velocity $u$, when the speed parameter is $v-u$.

To motivate our more general setting, consider the following example of a multi-agent distributed synchronization model.
We have $n$ agents (particles), each trying to advance its own state (particle location), but there are dependencies 
between the agents. 
(One specific motivation is a parallel simulation, where an agent is a processor, simulating 
an element, or a subsystem, of a larger physical system. An agent state is the local simulation time of the processor.
Cf. \cite{GSS96}.) 
A job arrival corresponds to selection of a group of agents, of size $d$, that are currently interdependent.
Upon such selection, the agents in the group, first, separately advance their states (locations) $W_1, \ldots, W_d$ forward by random distances (potential jump
sizes) $\xi_1,\ldots, \xi_{d}$, but then have to synchronize -- each agent $i$ rolls back its new state from $W_i +\xi_i$ to 
$\max\{W_i, \min_\ell (W_\ell + \xi_\ell)\}$.
(In terms of parallel simulation,  $\xi_i$ is the time duration by which processor $i$ can advance
its local simulation without affecting simulations run by other processors in the group.) 
If we assume that jobs arrive as a Poisson process, and the groups are selected uniformly at random,
we obtain a single-class unregulated (free) system with $(d,1)$-c.o.c. redundancy. 

Results of this paper, in particular, characterize the limits
$\lim_{n\to\infty} v_n$ of the average  
velocity $v_n$
for general free systems. This characterization is particularly sharp in the special case 
when the component sizes $(\xi^{(j)}_1,\ldots, \xi^{(j)}_{d_j})$ for each class are independent identically distributed,
with the distribution having increasing hazard rate.

\subsection{Paper contributions}

The main contributions of this paper are that we: (a) generalize and unify results for the left-regulated system, as compared to \cite{St2021-coc}; 
(b) obtain results for the free and left-regulated systems; and (c) our analysis demonstrates and exploits the connections between the systems with different types of regulation. Specifically, some of our main results, stated here informally, are as follows. (We give references to formal results, stated 
in Section~\ref{sec-main-results-new}.)
The results hold under the assumption of finite second moment of job component sizes
and very mild non-degeneracy assumptions.

\begin{itemize} 

\item For the free system, a unique ML fixed point exists for any speed $v$ within some non-empty range $[v_{min},v_{max}]$.
(Theorem~\ref{thm-free-wave-spectrum}(i))

\item For any speed $v < v_{min}$ 
there exists a unique ML fixed point for the right-regulated system
(Theorem~\ref{thm-free-wave-spectrum}(iii))

\item For 
any speed $v > v_{max}$ there exists a unique ML fixed point for the left-regulated system.
(Theorem~\ref{thm-free-wave-to-regulated-unique}(i))

\item Property $v_{min}=v_{max}$ holds under the additional condition of $D$-monotonicity (Definition~\ref{def-d-monotone}) on the component sizes distribution (Theorem~\ref{thm-free-wave-unique-cond}(ii)).
In turn, $D$-monotonicity holds when for each job class $j$ the component sizes are i.i.d. with a distribution having increasing hazard rate
(Definition~\ref{def-ihr}).

\item In the left-regulated system, SSAI (Definition~\ref{def-ssai})
holds for any speed $v > v_{max}$ (Theorem~\ref{thm-free-wave-to-ssai}). As $v \downarrow v_{\max}$,
the ML fixed point (on which stationary distributions concentrate) is such that the system load increases and converges to some $\bar \rho \le 1$.
If $\bar \rho < 1$ the fixed point converges to the minimum ML fixed point for the left-regulated system with speed $v_{max}$; the latter has load $\bar \rho$.
If $\bar \rho = 1$, the ML fixed point converges to the ML fixed point for the free system with speed $v_{max}$. An additional condition
(for which $D$-monotonicity is sufficient)
 is given, under which
$\bar \rho = 1$. Theorem~\ref{thm-free-wave-to-ssai} unifies and generalizes the main results of \cite{St2021-coc} on the SSAI
in left-regulated systems, as shown in Corollaries~\ref{thm-free-wave-to-ssai-cor1} and \ref{thm-free-wave-to-ssai-cor2}.

\item In the right-regulated system, SSAI (Definition~\ref{def-ssai})
holds for any speed $v < v_{min}$ (Theorem~\ref{thm-ssai-right}). As $v \uparrow v_{\min}$,
the ML fixed point (on which stationary distributions concentrate) converges to the ML fixed point 
for the free system with speed $v_{min}$. 

\item Let $v_n$ denote the average velocity at which the entire free $n$-particle system advances, in the model with zero speed parameter.
Then $\liminf_n v_n \ge v_{min}$ and $\limsup_n v_n \le v_{max}$ (Theorem~\ref{thm-vn-limit}). In the special case 
when for each job class $j$ the component sizes are i.i.d. with a distribution having increasing hazard rate
(Definition~\ref{def-ihr}), but with other assumptions being even milder than in Theorem~\ref{thm-vn-limit}, 
we have $\lim_n v_n = v_{min}=v_{max}$ (Theorem~\ref{thm-vn-iid-ihr}).

\end{itemize}

\subsection{Prior work}

There is a large amount of prior and ongoing work on redundancy models in queueing. 
(We note that in most of the prior work what we call job components are called job ``replicas,'' or ``copies.'')
C.o.c. is not the only form of redundancy and mean-field asymptotic regime is not the only regime of interest. 
See, for example, \cite{6025217e801240d598341bfdd12e8c43} for a general review of redundancy models in multi-server queueing systems.
A review of prior work more closely related to the present paper, i.e. to c.o.c. redundancy and/or mean-field asymptotics,
 can be found in \cite{St2021-coc,SnSt2020}, which in particular includes papers
 \cite{Vulimiri-2013,shah2015redundant,gardner2015reducing,Harcol-Balter-2017,adan2018fcfs,ayesta2018unifying,ayesta2019redundancy,VDK96,BLP2012-jsq-asymp-indep,St2021-coc,SnSt2020}.
 
 We note that prior work {\em proving} SSAI for redundancy systems in the mean-field regime is rather scarce. 
 It includes papers \cite{VDK96,BLP2012-jsq-asymp-indep} which prove results on SSAI under the extensively studied {\em power-of-$d$-choices} load balancing (component placement) scheme. Paper \cite{SnSt2020} proves SSAI for a system with multi-component jobs, which employs {\em water-filling} and/or {\em least-load}
 component placement mechanisms. Least-load mechanism includes the so-called {\em cancel-on-start} redundancy,
 which in turn includes {\em power-of-$d$-choices}. All systems studied in \cite{VDK96,BLP2012-jsq-asymp-indep,SnSt2020} are work-conserving.
 The water-filling mechanism includes c.o.c. redundancy, {\em but only in the special case of i.i.d. exponentially distributed component sizes.}
 Paper \cite{St2021-coc}, as already stated, studies the special case -- left-regulated system -- of the more general model in present paper and,
 in particular, obtains SSAI results -- this will be extensively discussed later.
 
While results {\em proving} SSAI for redundancy models are scarce, 
the {\em SSAI conjecture} is often employed to compute the corresponding ML fixed point (on which the stationary distribution concentrates as $n\to\infty$,
assuming SSAI holds), which is then used to
obtain estimates of the system performance metrics when the number of servers is large (cf. \cite{Vulimiri-2013,Harcol-Balter-2017,hellemans2019performance}). Therefore, for systems, for which we prove SSAI in this paper, we rigorously substantiate 
such computations and estimates. 

Recall that the unregulated (free) version of our model is the system of particles evolving on the real line, without any boundaries.
There is an extensive literature on the mean-field asymptotics for such systems, including mean-field limits, their traveling wave solutions (ML fixed points in our terminology), and SSAI. Cf. reviews in \cite{St2022-wave-lst, BaSt2023}.

\subsection{Paper layout}
\label{sec-layout} 

The rest of the paper is organized as follows. Basic notation and terminology, used throughout the paper, are introduced in Section~\ref{sec-basic-notation}. In Section~\ref{sec-formal-model} the model and the asymptotic regime are formally defined. 
In Section~\ref{sec-tools} we describe basic monotonicity and continuity properties of the model, as well as the notion of a reduced system,
which are used throughout the paper. In Section~\ref{sec-ml-new} we define mean-field limits and present those of their properties that generalize the
corresponding properties derived \cite{St2021-coc}. In Section~\ref{sec-res-coc-gen} we present preliminary definitions and results, generalizing those in \cite{St2021-coc}, which are related to SSAI. Section~\ref{sec-main-results-new} contains the main results of this paper. Section~\ref{sec-proofs-new} contains
proofs deferred from earlier sections. Concluding remarks are given in Section~\ref{sec-discussion}.

\section{Basic notation and terminology}
\label{sec-basic-notation}

We denote by $\R$ and $\R_+$ the sets of real and real non-negative numbers, respectively, and by $\R^n$ and $\R_+^n$ the corresponding $n$-dimensional product sets. By $\bar \R \doteq \R \cup \{\infty\} \cup \{-\infty\}$ we denote the two-point compactification of $\R$, where 
$\infty$ and $-\infty$ are the points at infinity and minus infinity, with the natural topology and a consistent with it metric,
so that $\bar \R$ [and $\R$] is complete and separable.
Analogously,
$\bar \R_+ \doteq \R_+ \cup \{\infty\}$ is the one-point compactification of $\R_+$. 
For a vector $q=(q_i) \in \R^n$,  $|q| \doteq \sum_i |q_i|$.
Inequalities applied to vectors [resp.  functions] are understood component-wise [resp. for every value of the function argument]. 
We say that a function is RCLL [resp. LCRL] 
if it is {\em right-continuous with left-limits} [resp. {\em left-continuous with right-limits}]. 
A scalar function $f(w)$ is called $c$-Lipschitz if it is Lipschitz continuous with constant $c$.
The sup-norm of a scalar function $f(w)$ is denoted $\|f(\cdot)\| \doteq \sup_w |f(w)|$; the corresponding convergence is denoted by $\stackrel{\|\cdot\|}{\longrightarrow}$. {\em U.o.c.} convergence means {\em uniform on compact sets} convergence, and is denoted by $\stackrel{u.o.c.}{\longrightarrow}$. We use notation: $a\vee b \doteq \max\{a,b\}$, 
$a\wedge b \doteq \min\{a,b\}$.
Abbreviation {\em w.r.t.} means {\em with respect to};
{\em a.e.} means {\em almost everywhere w.r.t. Lebesgue measure};
WLOG means {\em without loss of generality}; RHS and LHS means {\em right-hand side} and {\em left-hand side}, 
respectively.

We denote by $\tilde \cx$ [resp., $\cx$] the set of non-increasing RCLL functions $x= (x_w, ~w\in \R)$, [resp., $x=(x_w, ~w\in \R_+)$,]
taking values in $[0,1]$. An element $x \in \tilde \cx$ [resp., $x\in \cx$] is naturally interpreted as complementary distribution function
on $\bar \R$ [resp., $\bar \R_+$], with $1-x_w$ being the measure of $[-\infty,w]$ [resp., $[0,w]$]. An element $x \in \tilde \cx$ is {\em proper} if $x_\infty \doteq \lim_{w\to \infty} x_w =0$ and $x_{-\infty} \doteq \lim_{w\to -\infty} x_w =1$, and {\em improper} otherwise; in any case, $x_\infty$ and $1 - x_{-\infty}$ give the measures of $\{+\infty\}$ and $\{-\infty\}$, respectively. Analogously, an element 
$x \in \cx$ is proper if $x_\infty =0$ and  and improper otherwise. 
(With some abuse of standard terminology, in this paper we sometimes refer to elements of  $\tilde \cx$ and $\cx$ as distributions.)
Denote by $\tilde \cx^{pr}$ [resp., $\cx^{pr}$] the subset of proper elements of $\tilde \cx$ [resp., $\cx$].

We will equip the space $\tilde \cx$ with the topology of weak convergence of the corresponding distributions on $\bar \R$;
equivalently, $y \to x$ if and only if $y_w \to x_w$ for each $w \in (-\infty,\infty)$ where $x$ is continuous.
\sasha{Furthermore, on $\tilde \cx$ we consider the following metric $L$, consistent with the weak convergence topology.
We map an element $x \in \tilde \cx$ into the (proper or non-proper) distribution function $\phi=\phi(x)$ on $\R$, where
$\phi_s = 1- x_{-\log (1-s)}$ for $0 \le s < 1$, $\phi_s = 1- x_{\log (1+s)}$  for $-1 \le s < 0$, $\phi_s = 0$ for $s<-1$, $\phi_s = 1$ for $s\ge 1$.
Then $L(x,y) \doteq \hat L(\phi(x),\phi(y))$, where $\hat L$ is the Levy-Prohorov metric (cf.  \cite{Ethier_Kurtz}). 
It is easy to see that the convergence in metric $L$ is indeed equivalent to the weak convergence.}
Clearly, $\tilde \cx$ is compact.
Space $\cx$ inherits the topology and metric of space $\tilde \cx$ -- for these purposes we view any $x \in \cx$ as 
an element of $\tilde \cx$ via the convention that $x_w = 1$ for $w<0$. Clearly, $\cx$ is also compact.

For a function $x \in \tilde \cx$ [resp., $x\in \cx$] we denote by $x^{-1}=(x_u^{-1}, ~u \in [0,1]),$
its inverse, defined as
$$
x_u^{-1} = \sup \{w~|~x_w > u\},
$$
with the convention that $x_u^{-1} = -\infty$ [resp., $x_u^{-1} = 0$]
if the set in $\sup$ is empty. 

Unless explicitly specified otherwise, we use the following conventions regarding random elements and random processes.
A measurable space is considered equipped with a Borel $\sigma$-algebra, induced by the 
metric which is clear from the context. A random process $Y(t), ~t\ge 0,$ always takes values in a complete separable metric space (clear
from the context), and has RCLL sample paths; the sample paths are elements of the Skorohod $J_1$-space with the corresponding metric.

For a random process $Y(t), ~t\ge 0,$ we denote by $Y(\infty)$ the random value of $Y(t)$ in a stationary regime (which will be clear from the context). Symbol $\Rightarrow$ signifies convergence of random elements in distribution; $\stackrel{P}{\longrightarrow}$ means convergence in probability.
 {\em W.p.1} means {\em with probability one.}
{\em I.i.d.} means {\em independent identically distributed.}
Indicator of event or condition $B$ is denoted by $\bI(B)$. If $X,Y$ are random elements taking values in a set on which a partial order $\le$ is defined, then the stochastic order
$X \le_{st} Y$ means that $X$ and $Y$ can be coupled (constructed on a common probability space) so that $X \le Y$ w.p.1.

For an element $y\in \tilde \cx$, denote by $\bar y$ the mean of the corresponding distribution,
$$
\bar y = \int_{-\infty}^\infty u d[-y_u] = \int_0^\infty y_u du - \int_{-\infty}^0 (1-y_u) du,
$$
with the usual convention that $\bar y$ is well-defined when at least one of the 
integrals in the RHS is finite. 
When $\bar y$ is finite, denote by $\mathring y = (\mathring y_u,~u\in \R)$ the centered version of $y$, namely
$$
\mathring y_u = y_{u+\bar y}, ~u\in \R,
$$
and denote
$$
\Phi_\ell(y) = \int_{-\infty}^\infty |u|^\ell d[- \mathring y_u] = \int_{-\infty}^\infty |u-\bar y|^\ell d[-y_u], ~~\ell \ge 1,
$$
$$
\tilde \cx_\ell = \{y\in \tilde \cx ~|~ \mbox{$\bar y$ is finite}, ~ \Phi_\ell(y) < \infty\}, ~~\ell \ge 1,
$$
$$
\mathring{\cx}_1 = \{y\in \tilde \cx_1 ~|~ \bar y = 0\}.
$$

\section{Model}
\label{sec-formal-model}

\subsection{Multi-server queueing system with cancel-on-completion redundancy}

Consider a queueing model with $n$ identical servers, each processing its work at rate $1$.
The {\em workload} of a server at a given time $t$ is its unfinished work,
namely the time duration until the server becomes idle {\em assuming no new job arrivals after $t$.} 
There is a finite set $\cj$
of job classes $j$. 
Jobs  arrive as Poisson process of rate $\lambda^n n$, for some $\lambda^n \ge 0$.
The probability that an arriving job is of class $j$ is $\pi_j^n \ge 0$, $\sum_j \pi_j^n = 1$.
We will also use notation $\sigma_j^n = \lambda^n \pi_j^n$ for the class $j$ arrival rate per server.

Each job class $j$ has three parameters: integers $k_j$ and $d_j$ %(they will not depend on $n$)
such that $1\le k_j \le d_j$, and an exchangeable probability distribution $F_j$ on $\R_+^{d_j}$. 
(Exchangeability of $F_j$ means that it is invariant w.r.t. permutations of components.)
When a class $j$ job arrives, $d_j$ servers are selected uniformly at random (without replacement); these servers form the {\em selection set} of the job. The job places $d_j$ {\em components} on the selected servers; the component {\em sizes} (potential amounts of required service) 
$(\xi^{(j)}_1,\ldots, \xi^{(j)}_{d_j})$ are drawn according to distribution $F_j$, independently of the process history up to the job arrival time.
Each server processes its work (components of different jobs) in the First-Come-First-Served (FCFS) order. 
A job of class $j$, to be completed, requires $k_j$ (out of $d_j$) components to be processed, and as soon as $k_j$ components of the job {\em complete} their service (i.e. receive the amounts of service equal to their sizes), the unfinished service of the remaining $d_j -k_j$ components of the job is ``canceled'' and 
the job immediately leaves the system;
we will call this {\em $(d_j,k_j)$-c.o.c. redundancy}.
Throughout the paper we use notation $\bar d \doteq \max_j d_j$.

Clearly, when a job arrives and selects a server, the {\em actual} amount $\kappa$ of workload added to the server may be {\em smaller} 
than the size $\xi$ of the job component placed on the server, i.e. $0 \le \kappa \le \xi$,
due to possible cancellation (of a part or all) of the component service.
Moreover, $\kappa$ depends not only on the realization of the job component sizes, 
but also on the current workloads of the selected servers or, more specifically, on their {\em relative} workloads, i.e. workload differences.
(Also note that due to FCFS discipline at each server, the actual amounts of workload added to selected servers do {\em not} depend on the future job arrivals; so that new workloads are well-defined after each new arrival.)
Because of this, c.o.c. redundancy in general has a property, which substantially  complicates its analysis, namely 
it is {\em non-work-conserving} in the following sense.
Suppose a class $j$ job arrives. Then, in general, the average total amount of new workload added by the arriving job to 
the selected servers
depends not only on the distribution $F_j$ of the component sizes,
but also on the relative workloads of the selected servers.
In particular, the average rate per server at which new workload is added to the system is not known a priori.
(If for each class $j$ either $k_j=d_j$ or the component sizes are i.i.d. exponential, the work-conservation does hold --
the average total amount of new workload added by an arriving class $j$ job depends on $j$ only; see \cite{St2021-coc,SnSt2020}.)

In this paper we will primarily use the following ``particle'' language to describe the system dynamics. Namely, we will identify each server with a ``particle,'' the server workload with ``particle location,''  and the workload evolution with ``particle movement.''
The basic dynamics of a particle (server workload) is as follows. If/when a job arrival adds workload $\kappa \ge 0$ to the server,
the particle jumps ``right'' by $\kappa$. Between job arrivals, the particle (workload) moves ``left'' at the constant speed $1$ until and unless it ``hits'' the left ``regulation boundary'' at $0$; if the particle does hit the boundary $0$, it stays at it until and unless it jumps right due to a job arrival.

\subsection{More general class of systems}

For both applications and the analysis, it is convenient to view the c.o.c. redundancy model above within a much more general class of models, which we now define. The job classes $j$ and their parameters $k_j$, $d_j$, $F_j$, $\sigma^n_j$ have the same meaning. 
(Parameters $\lambda^n= \sum_j \sigma^n_j$ and $\pi^n_j = \sigma^n_j/\lambda^n$ are determined by $\{\sigma^n_j\}$.)
Additional features of a more general model are as follows:

\begin{itemize}
\item The constant speed at which each particle moves to the left (between jumps and unless it is on a left boundary)
 is not necessarily $1$, but some $v \in [0,\infty)$. This speed $v$ is an additional model parameter.

\item Instead of a fixed left boundary at $0$, the system may have: both left boundary at point $A$ and right boundary at point $B \ge A$;
or only one -- left or right -- boundary; or have no boundaries at all. 
If there is a right boundary at point $B$, this means that if, after following the c.o.c. procedure upon a job arrival, some selected particles would ``jump over'' point $B$, they land exactly at $B$ instead.
\end{itemize}

System with finite left boundary $A$ will be called {\em left-regulated}, and without left boundary (which can be viewed as $A= -\infty$) -- {\em left-unregulated.} The terms {\em right-regulated}, {\em right-unregulated}, {\em both-sides-regulated}, etc., will have analogous self-explanatory meanings. 
{\em Unregulated} system -- with neither left nor right boundary -- will also be called {\em free} (as in \cite{St2021-coc}).
The interval between the left and right boundaries $A$ and $B$ is sometimes called the system {\em frame} (as in \cite{St2021-coc}).
For example, frame $(-\infty,0]$ means left-unregulated system with right boundary at $0$; finite frame $[A,B]$ means both-sides-regulated system;
and so on.

The class of systems we described so far is such that each particle moves left at a constant speed $v$ (between jumps and unless on a left boundary)
and the boundaries (if any) are fixed. Such models we will call {\em anchored}.

For both applications and analysis, for an anchored system it is sometimes convenient to consider its {\em unanchored} version,
where particles do not move left between jumps, but instead the boundaries (if any) move right at the constant speed $v$. 
Clearly, anchored and the corresponding unanchored system are equivalent, in the sense that the former is the latter viewed in the system of coordinates moving right at constant speed $v$. (In other words, the unanchored system has a frame moving right at speed $v$, and the anchored version 
is obtained by considering the system evolution with respect to its frame.) Obviously, in the special case of $v=0$ the anchored and unanchored systems are indistinguishable. The case of $v>0$ will of primary interest in this paper, but some of the results apply more generally to $v \ge 0$.

The multi-agent synchronization model described in Section~\ref{sec-intro}, is an example of an unanchored unregulated (free) system. In such systems the speed parameter $v$ is irrelevant. As discussed earlier, for such a system one of the interesting characteristics is the average velocity $v_n$ at which the entire population of particles advances to the right. In this paper we will use the term ``velocity'' to refer to the average rate at which the particle population (or a mean-field limit, defined later) advances to the right; this is to make a distinction from the term ``speed,'' which reserve for the speed parameter $v$. Usually, we consider velocities in unanchored systems.

\subsection{Assumptions on job component size distributions}

Throughout the paper we always make the following Assumption~\ref{cond-finite-mean}, without stating it explicitly (although, for some results 
even this assumption can be relaxed.)

\begin{assumption}
\label{cond-finite-mean}
(i) For all $j$, 
\beql{eq-finite-mean}
\E \xi^{(j)}_1 < \infty ~~ \mbox{(and then $\E \sum_{i=1}^{d_j} \xi^{(j)}_i = d_j \E \xi^{(j)}_1 < \infty$),}
\eeql
and WLOG, $\E \xi^{(j)}_1 >0$.

(ii) There exists class $j$ with $F_j$ such that
\beql{eq-non-triv1}
\pr \{ \sum_i \bI(\xi^{(j)}_i =0) \ge k_j \} < 1.
\eeql
\end{assumption}

Assumption~ \ref{cond-finite-mean}(ii) implies that, uniformly on the system state at the time of a 
job arrival, the expected total amount of 
workload (the expected total displacement of particles) that the job brings to the system is at least some $\delta>0$.
If Assumption~ \ref{cond-finite-mean}(ii) does {\em not} hold, the system behavior can be degenerate in that when the system reaches the state with
all particles co-located, the particles will remain co-located from that time on.

Assumptions~\ref{cond-second-moment}, \ref{cond-ntriv-add-ii} and \ref{cond-ntriv-add-i} on the job component sizes distributions, which we list next, are {\em not} required for
many results, and will be stated explicitly (along with all other required assumptions) when needed.

\begin{assumption}
\label{cond-second-moment}
For any $j$, $\E(\xi_1^{(j)})^2 < \infty$ (and then $\E(\sum_i \xi_i^{(j)})^2 < \infty$).
\end{assumption}

\begin{assumption}
\label{cond-ntriv-add-ii}
There exists a class $j$, with $k_j < d_j$ and the joint distribution $F_j$ of the component sizes such that
$$
\pr \{ \sum_i \bI(\xi^{(j)}_i =0) \ge k_j \} < 1.
$$
\end{assumption}

 Assumption~\ref{cond-ntriv-add-ii} is a slightly stronger version of the non-degeneracy
 Assumption~\ref{cond-finite-mean}(ii), additionally requiring that $k_j < d_j$ for the class $j$, and 
 is very non-restrictive. As an example, it holds when there is a class $j$ with $k_j < d_j$ and i.i.d. component sizes.

\begin{assumption}
\label{cond-ntriv-add-i}
There exists a class $j$, with $k_j < d_j$ and the joint distribution $F_j$ of the component sizes such that
$$
\P\{ \max_i \xi^{(j)}_i - \xi^{(j)}_{i_1} >0\} > 0,
$$
 where $\xi^{(j)}_{i_1}$ is the $k_j$-th smallest among the component sizes $\xi^{(j)}_i$.
\end{assumption}

Assumption~\ref{cond-ntriv-add-i} also holds in most cases of interest. For example, it automatically holds when there 
is a class $j$ with $k_j < d_j$ and i.i.d. component sizes, with a component size distribution {\em not} concentrated on a single point.
The assumption does not hold, for example, when each class $j$ has all components of equal (possibly random) size.

 \subsection{Asymptotic regime. Mean-field scaled process} 
 \label{sec-asymp-regime}
 
 We consider a sequence of systems -- anchored, to be specific --  with $n\to\infty$, $\lambda^n \to \lambda \ge 0$
 and $\pi_j^n \to \pi_j$ for all $j$ (or, equivalently, $\sigma_j^n \to \sigma_j$ for all $j$),
 while all other system parameters remain fixed. 
 
From now on, the upper index $n$ of a variable/quantity will indicate that it pertains to the system with $n$ servers,
or $n$-th system.
Let $W_i^n(t)$ denote the location
particle $i$ at time $t$ in the $n$-th system. 
Clearly, for each $n$ the process $W^n(t) = (W^n_1(t), \ldots, W^n_n(t)), ~t\ge 0,$ is Markov with state space $\R^n$.
We adopt the convention that its sample paths are RCLL. 

Consider also the following {\em mean-field} scaled quantities:
\beql{eq-x-def2}
x^n_w(t) \doteq (1/n) \sum_i \bI\{W_i^n(t)> w\}, ~~ w \ge 0.
\end{equation}
That is, $x^n_w(t)$ is the fraction of particles $i$ with $W_i^n(t)> w$, i.e. those located to the right of point $w$.
Then
 $x^n(t)=(x^n_w(t), ~w\ge 0)$ describes the empirical distribution of particles' locations, it is a projection of $W^n(t)$, and $x^n(t), ~t\ge 0,$ is a Markov process.

Clearly, for any $n$ and $t$, $x^n(t) \in \cx^{pr}$. Denote by $\cx^{(n)} \subset \cx^{pr}$
the state space $\cx^{(n)}$ of the Markov process $x^n(t), ~t\ge 0.$ 
Therefore,
for any $n$, we can and will view
$x^n(t), ~t\ge 0$ as a Markov process with (common) state space 
$\cx^{pr} \subset \cx$,
and with sample paths being RCLL functions of $t\ge 0$ (taking values in $\cx^{pr}$).

We say that the process $W^n(\cdot)$ [resp. $x^n(\cdot)$] is {\em stable} if it is {\em positive Harris recurrent}.
For a general definition of positive Harris recurrence, cf. \cite{Bramson-book}.
Obviously, $W^n(\cdot)$ is stable if and only if $x^n(\cdot)$ is.
As far as positive recurrence is concerned, for the unregulated (free) system it only makes sense to consider 
the centered process $\mathring x^n(\cdot)$ instead of $x^n(\cdot)$.

In the special case of a regulated system, the process has simple structure -- it is irreducible, possessing a fixed state being a renewal atom.
For example, suppose the system is left-regulated at boundary point $0$. (It also may or may not have a right boundary.)
Then the state $(0,\ldots,0)$ of $W^n(\cdot)$, with all particles at the left boundary, is a renewal atom; 
this state we will call ``empty'', because in terms of the multi-server queueing system it corresponds to all servers being idle.
The corresponding renewal atom for the process $x^n(\cdot)$ is the state $x^\emptyset$, defined by
\beql{eq-empty-def}
x^\emptyset_w = \bI\{w < 0\};
\eeql
it will also be called the empty state. The irreducibility of $W^n(\cdot)$ [resp., $x^n(\cdot)$] follows from the facts that 
the empty state is reachable from any other, and that the job arrival process is Poisson.

For a regulated system the stability (positive Harris recurrence) simply means 
that the renewal atom is reachable from any other state w.p.1 and the 
expected time to return to the renewal (after leaving it) is finite.

Note that, trivially, the both-side-regulated (finite-frame) system is always stable. If the system is unregulated on at least one side (i.e., has infinite frame),
than due to the fact that the system is non-work-conserving, the stability conditions are not automatic and is non-trivial.

Stability of the process $x^n(\cdot)$ (or,  $\mathring x^n(\cdot)$ for the free system)
implies that it has unique stationary distribution. If the process is stable,
let $x^n(\infty)$ be a random element with values in $\tilde \cx^{pr} \subset \tilde \cx$,
whose distribution is the stationary distribution of the process; in other words,
it is a random process state in the stationary regime. 

Later in the paper we will need some additional notation, associated with the spaces $\tilde \cx$ and $\cx$.
The ``infinite'' state 
$x^{**} \in \cx$ is defined by
\beql{eq-inf-def}
x^{**}_w \equiv 1.
\eeql
It corresponds to all particles located at $+\infty$.

\subsection{Steady-state asymptotic independence definition}
\label{sec-ssai-def}

\begin{definition}
\label{def-ssai}
Consider our asymptotic regime with $n\to\infty$. 

(i) We say that {\em steady-state asymptotic independence} (SSAI) holds for speed $v>0$ for the anchored system, which is either left-regulated, right-regulated or both-side-regulated, if the Markov process $x^n(\cdot)$ is stable (positive recurrent) for all large $n$ and
\beql{eq-conv-main}
x^n(\infty) \Rightarrow x^{*},
\eeql
for some fixed $x^* \in \tilde \cx^{pr}$.

(ii) We say that {\em steady-state asymptotic independence} (SSAI) holds for the unregulated (free) system, if the centered Markov process $\mathring x^n(\cdot)$ is stable (positive recurrent) for all large $n$ and
\beql{eq-conv-main-free}
\mathring x^n(\infty) \Rightarrow x^{*},
\eeql
for some fixed $x^* \in \tilde \cx^{pr}$.
\end{definition}

For the anchored system with frame $[0,\infty)$ (i.e. left-regulated at $0$ and unregulated on the right), the following property will be called
{\em SSAI for the full range of loads} (SSAI-FRL): there exists $\underline{v}$ such that for any 
$v > \underline{v}$ the convergence $x^n(\infty) \Rightarrow x^{*,v}$ holds for some fixed $x^{*,v} \in \tilde \cx^{pr}$ and
$x^{*,v}_0 \uparrow 1$ as $v \downarrow \underline{v}$. The name SSAI-FRL is natural, because, in terms of the multi-server queueing system,
$x^{*,v}_0$ is the limiting fraction of occupied servers -- the system load. The notion of SSAI-FRL was introduced in \cite{St2021-coc}; see 
\cite[section 1.1]{St2021-coc} for a discussion of this property, which is not automatic -- and may be challenging to prove -- for non-work-concerving systems.

\section{Useful tools}
\label{sec-tools}

\subsection{Basic monotonicity and continuity}
\label{sec-basic-monot-cont}

Our model satisfies basic path-wise monotonicity and continuity properties, all stemming 
from the following simple fact. 

\begin{lem}[Lemma 1 in \cite{St2021-coc}]
\label{lem-monotone-cont-basic}
Consider a fixed set of particles, selected by a class-$j$ job, and labeled WLOG by $1, 2, \ldots, d_j$, with (deterministic) locations given by vector
$W=(W_1, \ldots, W_{d_j})$. Let $\xi^{(j)} = (\xi^{(j)}_1, \ldots, \xi^{(j)}_{d_j})$ be the (deterministic) vector of their potential jump sizes. Denote 
by $\hat W=(\hat W_1, \ldots, \hat W_{d_j})$ the vector of the particles' locations after the job arrival. Then $\hat W$ is monotone non-decreasing and continuous in both $W$ and $\xi^{(j)}$.
\end{lem}

This fact easily implies properties like: two versions of the same process, one starting from a smaller state than another, can be coupled 
so that the state of the former is smaller than the state of the latter at all times. This property obviously prevails, for example, if, in addition, the distribution of potential jump sizes in the former process is dominated by that in the latter, and/or if movement of particles in the former system is ``impeded'' by a fixed or moving right boundary, and/or if movement of particles in the latter system is ``helped'' by a fixed or moving left boundary, and so on. 
Similarly, for example, if two processes start from the same initial state and the distribution of potential jump sizes in the former converges to those in the latter, then the former process converges (in appropriate sense) to the latter; analogous continuity properties hold as well.
Furthermore, such monotonicity and continuity properties are naturally inherited by the mean-field limits, which we introduce in Section~\ref{sec-ml-new}.
All such properties are straightforward -- we will not spell out all of them. (Some are spelled out in section 6 of \cite{St2021-coc}.)
In the rest of this paper, we simply refer to them, where necessary,
as {\em monotonicity} and/or {\em continuity}. 

\subsection{Reduced systems and corresponding job subclasses}
\label{sec-reduced-and-subclasses}

We now define {\em reduced} systems and job {\em subclasses} -- the notions used in the analysis throughout the paper. One of the primary uses of the reduced system notion is as follows: it allows us to interpret a ``truncation'' of an ML fixed point (defined in Section~\ref{sec-ml-new}) as an ML fixed point of a reduced system, which in turn allows us to use properties of the latter to analyze the original ML fixed point.  (This is done, for example,
in  Theorem~\ref{thm-free-wave-unique-for-speed}, as well as in other places.)

Consider a given system, with fixed $n$, with the set of job classes $\{j\}$ and corresponding parameters $\{k_j, d_j, F_j, \sigma_j\}$ for each $j$.
(We will not use superscript $n$ in this subsection to avoid clogging notation.) Let two parameters, $\epsilon \ge 0$ and $\delta \ge 0$, with $\epsilon + \delta < 1$, be fixed. Suppose that fraction $\epsilon$ of all particles are located at $+\infty$ (and therefore alway remain there) and fraction $\delta$ 
of all particles are located at $-\infty$ (and therefore alway remain there). The remaining fraction $1-\epsilon-\delta$ of particles are ``regular'' -- they are initially (and then always) at located at points in $\R$. What is the evolution of the ``regular'' particles? It can be equivalently described
as the evolution of the following {\em $(\epsilon,\delta)$-reduced} system.
First of all, in the $(\epsilon,\delta)$-reduced system there are $(1-\epsilon-\delta)n$ particles, instead of $n$. 
For each class $j$ in the original system, there are multiple classes in the reduced system. 
Specifically, consider two fixed integers $\ell$ and $m$, such that $0 \le \ell < k_j$ and $\ell < m \le d_j$. 
If a class $j$ job arrives into the original system, {\em and it happens that it selects exactly $\ell$ ``$-\infty$-particles'' and exactly $d_j-m$ ``$+\infty$-particles,''} then, in terms of the reduced system, this is equivalent to the arrival of a job of the class $j^{m,\ell}$ with the following parameters
(we will use notation $\omega =j^{m,\ell}$ for brevity): $k_\omega = k_j -\ell$, $d_\omega= m-\ell$; the distribution $F_\omega$ of the component sizes 
is the projection of $F_j$ on the distribution of first $m-\ell$ components; the arrival rate (per server) is 
\beql{eq-sig-omega}
\sigma_\omega = \sigma_j \pr\{ \mbox{a job $j$ selects exactly $\ell$ ``$-\infty$-particles'' and exactly $d_j-m$ ``$+\infty$-particles''} \} / (1-\epsilon-\delta).
\eeql
Indeed, from the point of view of reduced system, the selected $d_j - m$ ``$+\infty$-particles'' can be ignored, and the selected  $\ell$ ``$-\infty$-particles'' reduce the number of particles whose movement needs to be completed from $k_j$ to $k_\omega=k_j - \ell$; thus, in the reduced system we only need to consider the selection of $m-\ell$ particles, of which the movement of $k_j-\ell$ needs completion. The renormalization by factor 
$(1-\epsilon-\delta)^{-1}$ in \eqn{eq-sig-omega} is needed because the number of particles is $(1-\epsilon-\delta)n$ instead of $n$.
The constraint $\ell < k_j$ is because if $k_j$ or more ``$-\infty$-particles'' are selected by a class $j$ in the original system,
this has no impact on the reduced system; similarly, the constraint $\ell < m$ is because, otherwise, no particle in the reduced system is selected.

Classes $j^{m,\ell}$ corresponding to the original class $j$, we will call {\em subclasses} of $j$. Obviously, the subclass $j^{d_j,0}$
is equal to the class $j$ itself. 

We note that the job classes in the $(\epsilon,0)$-reduced system are all subclasses of the form $j^{m,0}$, for all $j$, in the original system.
(In \cite{St2021-coc} only $(\epsilon,0)$-reduced system and corresponding subclasses $j^{m,0}$ are considered.)
In the $(0,\delta)$-reduced system the classes are subclasses of the form $j^{d_j,\ell}$, for all $j$, in the original system.

A reduced system is a well-defined system within our model. Therefore, all notions, constructions and results still apply.
In particular, regulated and unregulated versions can be considered. Also, 
for a reduced system we can consider the same asymptotic regime, with  $n\to\infty$, %and with some of the parameters converging rather being fixed, 
as for the original one. Note that
\beql{eq-sig-omega-lim}
\lim_{n\to\infty} \frac{\sigma_\omega}{\sigma_j} =  \frac{d_j !} {\ell ! (m-\ell)! (d_j-m) !} \delta^\ell (1-\epsilon-\delta)^{m-\ell} \epsilon^{d_j-m} 
(1-\epsilon-\delta)^{-1}.
\eeql
In particular, all notions and results associated with the asymptotic regime, such as mean-field limits, described later, apply to reduced systems as well.

\section{Mean-field limits}
\label{sec-ml-new}

In this section we define the process mean-field limits (ML) and some of their properties. These definitions and properties generalize those in \cite{St2021-coc}.
We only provide proofs which are substantially different from the corresponding proofs in \cite{St2021-coc}. Further ML properties, which are some of the main results of this paper, are derived in Section~\ref{sec-main-results-new}.

For any element $x \in \tilde \cx$, let us use a slightly abusive notation $x_{(-\infty,w]}$ for  
 $(x_u, ~u\in (-\infty,w])$, i.e. for function $x_\cdot$ with the domain truncated to $(-\infty,w]$. 
 (Here an element $x \in  \cx$ is treated as an element of $\tilde \cx$, via convention $x_w =1$ for $w <0$.)
 
 Suppose we randomly choose a job class
 $j$ according to the distribution $\{\pi_j\}$, and then assume a class $j$ job selects $d_j$ particles, whose locations are chosen randomly and independently according to distribution $x$. Then, let
 $$
 h(x_{(-\infty,w]})
 $$
 denote the expected number of those particles which, as a result, will cross level $w$ from left to right
 or, more precisely, will enter the interval $(w,\infty]$. 
Clearly, the values of $x_u$ for $u>w$ have no effect on this expectation -- that is why the expectation depends only on $x_{(-\infty,w]}$.

There is an alternative, equivalent way to define functional $h$, which will also be useful.
Let us use notations: $\alpha = \lambda \sum_j \pi_j d_j = \sum_j \sigma_j d_j$.  
For a fixed $x \in \tilde \cx$, $w\in \R$, and a job class $j$, consider CDF $J_j^{(x,w)}(u), ~u\ge 0,$ defined as follows.
Consider an (abstract) ``tagged'' particle, located at $w$, selected by a class $j$ job; suppose the remaining $d_j -1$ particles, selected 
by this job, have random independent locations chosen according to distribution $x$;
then $J_j^{(x,w)}(\cdot)$ is the CDF of the random actual jump size of the tagged particle. 
Denote $\bar J_j^{(x,w)}(u) = 1-J_j^{(x,w)}(u)$.
Then, it is easy to check that
\beql{eq-h-represent}
\lambda h(x_{(-\infty,w]}) = \alpha \sum_j \frac{\pi_j d_j}{\sum_\ell \pi_\ell d_\ell} \int_{-\infty}^w [-dx_u] \bar J_j^{(x,u)}(w-u)
= \sum_j \sigma_j d_j \int_{-\infty}^w [-dx_u] \bar J_j^{(x,u)}(w-u),
\eeql
or, using notation 
$$
J^{(x,w)}(u) = \sum_j \frac{\pi_j d_j}{\sum_\ell \pi_\ell d_\ell} J_j^{(x,w)}(u), ~~~ \bar J^{(x,w)}(u) = 1-J^{(x,w)}(u),
$$
$$
\lambda h(x_{(-\infty,w]}) = \alpha \int_{-\infty}^w [-dx_u] \bar J^{(x,u)}(w-u).
$$
Obviously (see the RHS of \eqn{eq-h-represent}), $\lambda h(x_{(-\infty,w]})$ is monotone non-decreasing w.r.t.
$\{\sigma_j\}$, where, recall, $\lambda = \sum_j \sigma_j$.
By basic monotonicity and coupling, if $x_{(-\infty,w]} \le y_{(-\infty,w]}$ and $x_w = y_w$, then $h(x_{(-\infty,w]})  \le h(y_{(-\infty,w]})$.
Finally, by the same argument as in lemma 12 in \cite{St2021-coc},  uniformly in $w$,
the functional $h(x_{(-\infty,w]})$ on $\tilde \cx$ is $\bar d^2$-Lipschitz ($\bar d = \max_j d_j$), w.r.t. the sup-norm on 
 $x_{(-\infty,w]}$, namely
 \beql{eq-h-Lipschitz}
 | h(x_{(-\infty,w]}) - h(y_{(-\infty,w]})| \le \bar d^2 \sup_{u \le w} | x_u - y_u |, ~~\forall x,y \in \tilde \cx.
 \eeql

\begin{definition}
\label{def-mm-new}
For an unanchored unregulated (free)  system, a function $x_w(t), ~w\in \R, ~t\in \R_+,$ will be called a {\em mean-field model} (MM), if it satisfies the following conditions.\\
 (a) For any $t$, $x(t) = (x_w(t), w\in \R) \in \tilde \cx^{pr}$.\\
(b) For any $w$, $x_w(t)$ is non-decreasing $c$-Lipschitz in $t$, with constant $c$ dependent only on the system parameters.\\
(c) For any $w$, for any $t$ where the partial derivative
$(\partial/\partial t) x_w(t)$ exists (which is almost all $t$ w.r.t. Lebesgue measure, by the Lipschitz property), equation 
\beql{eq-dyn-trans-new}
\frac{\partial}{\partial t} x_w(t) = \lambda h(x_{(-\infty,w]}(t))
\eeql
holds. 

For an unanchored right-regulated system with speed $v>0$ (i.e., the right boundary $B$ moving as $B=B_0+vt$), the MM definition is same as for the 
free system, except $x_w(t)=0$ for $w \ge B = B_0+vt$, and 
(c) 
holds for $w < B = B_0+vt$.

For an unanchored left-regulated system with speed $v>0$ (i.e., the left boundary $A$ moving as $A=A_0+vt$), the MM definition is same as for the 
free system, except $x_w(t)=1$ for $w < A = A_0+vt$, and 
(b) and (c) hold for $w > A = A_0+vt$.
\end{definition}

\begin{prop}
\label{prop-ml-result}
For an unanchored system with speed $v$ (which may or may not be right- or left-regulated)
we have 
the following.
Assume that the initial conditions are such that $x^{n}(0) \to x(0)$
for some fixed $x(0) \in \tilde \cx$.  
Then, there exists a deterministic trajectory $x(t), ~t\ge 0,$ with values in $\tilde \cx$ and initial state $x(0)$, such that
\beql{eq-conv-to-ml}
(x^n(t), ~t\ge 0) ~ \stackrel{P}{\longrightarrow} (x(t), ~t\ge 0).
\eeql
(That is, the sequence of processes $x^{n}(\cdot)$, with trajectories in the Skorohod space of functions taking values in $\tilde \cx$, 
converges 
in probability to the deterministic trajectory $x(\cdot)$.)
As a function of $t$ (i.e., as an element of Skorohod space), $x(t)$ is continuous. 
This trajectory $x$, which will be called the {\em mean-field limit} (ML), is the unique MM for this initial state $x(0)$.
\end{prop}

Proof of Proposition~\ref{prop-ml-result} is in Section~\ref{sec-ml-proof}.

\begin{definition}
\label{def-ml-any-speed}
For an anchored system with speed $v > 0$ (which may or may not be right- or left-regulated), the trajectory $x(t)=(x_w(t))$ will be called the ML
if the trajectory $y(t)=(y_w(t))$, with $y_w(t)=x_{w-vt}(t)$,
is the ML for the corresponding unanchored system with speed $v$.
\end{definition}

\begin{definition}
\label{def-ml-fp}
For an anchored system with speed $v > 0$,
function $x=(x_w)$ is called ML-FP if $x(t)=x$ is an ML. 
\end{definition}

If $x$ is an ML-FP for the anchored left-regulated system with frame $[A,\infty)$, the value $x_A$ at the left boundary will be called the {\em load} of $x$
(because this value corresponds to the fraction of occupied server in the multi-server queueing system).

\begin{definition}
\label{def-de-fp}
For an anchored system with speed $v > 0$, with left and right boundaries (if any) at points $A$ and $B$,
function $x=(x_w)$ is called DE-FP  if it satisfies the following conditions:\\
 (a) $x \in \tilde \cx^{pr}$, with $x_w = 1$ for $w<A$ and $x_w=0$ for $w\ge B$.                     \\
(b) $x_w$ is non-increasing $C$-Lipschitz, with constant $C$ dependent only on the system parameters, in the interval $[A,\infty)$.\\
(c) For almost all $w\in (A,B)$ w.r.t. Lebesgue measure,  
we have
\beql{eq-de-fp-new}
-v \frac{d}{dw} x_w = \lambda h(x_{(-\infty,w]}).
\eeql
\end{definition}

From Proposition~\ref{prop-ml-result} and the above definitions we obtain 

\begin{cor}
\label{prop-de-fp}
For an anchored system with speed $v > 0$,
function $x$ is an ML-FP if and only if it is DE-FP.
\end{cor}

Obviously, if $x$ is an ML-FP (=DE-FP) for an anchored system with speed $v>0$, then $x_{w-vt}$ is an ML in the corresponding unanchored system,
and it is a traveling wave, propagating to the right with velocity $v$. 
 For this reason, an ML-FP (= DE-FP) can also be called a {\em traveling wave shape (TWS).} 

Note that, if $x$ is an ML-FP in the anchored free system with speed $v>0$, then from \eqn{eq-de-fp-new} we obtain
(compare to (9) in \cite{St2021-coc})
\beql{eq-speed-equal-workload-arrivals}
v/\lambda =  \int_{-\infty}^{\infty} h(x_{(-\infty,u]}) du,
\eeql
where the RHS is the expected total displacement of selected particles upon a job arrival such that
the selected particles' locations are i.i.d. with distribution $x$.

For a given anchored system, 
let us define the following operator $\cal A$ which maps $\tilde \cx$ into itself. 
Let a distribution (``environment'') $x \in \tilde \cx$ be fixed; $x$ may or may not be proper. 
Consider a single particle moving in the frame $[A,B]$ (where the left boundary $A$ may be at $-\infty$, and right boundary $B$ may be at $+\infty$)
as follows. Unless and until the particle jumps, it moves left at speed $v$; if/when the particle ``hits'' the left boundary $A$, it stays there until next jump; there is a rate $\alpha$ Poisson process of time points when the the particle jumps forward; if the jump occurs when the particle location is $w$, the jump size distribution is $J^{(x,w)}(\cdot)$, independent of the process history (besides current location $w$); however, if the particle would jump over point $B$,
it lands exactly at $B$ instead. Then ${\cal A} x \in \tilde \cx$ is the stationary distribution of the particle. Note that, if a particle is located at $-\infty$ or $+\infty$, it stays there forever. Consequently,  $({\cal A}x)_\infty = x_\infty$ and $({\cal A}x)_{-\infty} = x_{-\infty}$.

\begin{definition}
A fixed point $x$ of the operator $\cal A$ we will call an {\em OP-FP}.
\end{definition}
   
  \begin{lem}
  \label{lem-fsp-is-op-fp}
If $x$ is a ML-FP, then it is an OP-FP. 
  \end{lem}
  
  {\em Proof}  is essentially same as that of lemma 10 in \cite{St2021-coc}.

  \begin{lem}
  \label{lem-op-fp-iff-de-fp}
   Any OP-FP $x$ is a DE-FP. 
  \end{lem}

   {\em Proof}  is essentially same as that of lemma 11 in \cite{St2021-coc}.

      \begin{cor}
  \label{cor-fp}
 An element $x \in \tilde \cx$ is a ML-FP if and only if it is an OP-FP, and if and only if it is a DE-FP.
  \end{cor}

\begin{lem}
  \label{lem-de-sol-unique}
 For any fixed initial condition $x_{(-\infty,0]}$, where $x_0\in [0,1]$ and the cases $x_{0-}=1$ and $x_{-\infty} < 1$ allowed,
 the solution $x_w, w>0,$ of the functional differential equation \eqn{eq-de-fp-new} exists, is unique and its dependence 
 on $x_{(-\infty,0]}$ is continuous. (Here the solution is in the sense that \eqn{eq-de-fp-new} holds for almost all $w>0$
 w.r.t. Lebesgue measure, and we consider the solution up to the point $w^*$ where it ``hits'' $w$-axis, if it does at all.)
\end{lem}

Lemma~\ref{lem-de-sol-unique} basically says that a DE-FP is uniquely determined by its ``left tail.''
 Its proof is essentially same as that of lemma 13 in \cite{St2021-coc}. 

The following general continuity ML-FP continuity property easily follows from the previous definitions and results  
in this section. 

\begin{lem}
  \label{lem-cont}
  Consider ML-FPs in anchored systems.
  Suppose we have a sequence, indexed by $\ell \to \infty$, of an anchored system parameters 
  $v^{(\ell)}>0$, $\{\sigma_j^{(\ell)}\}$, $\{F_j^{(\ell)}\}$, $A^{(\ell)}$, $B^{(\ell)}$,
  converging, as $\ell\to\infty$, to the limiting set of parameters $v>0$, $\{\sigma_j\}$, $\{F_j\}$, $A$, $B$, with $-\infty \le A  < B \le \infty$. 
  (Index $\ell$ is {\em not} the number $n$ of particles in a pre-limit system. Instead, it indexes ML-FPs for different parameter sets.
  Also, $A^{(\ell)}$, $B^{(\ell)}$ may be finite or infinite.)
  Suppose, for each $\ell$ there exists an ML-FP $x^{(\ell)} \in \tilde \cx$ (proper or not proper on either side) for the corresponding parameter set,
  and $x^{(\ell)} \to x \in \tilde \cx$. Then $x$ is an ML-FP for the limiting set of parameters.
\end{lem}

\begin{lem}
  \label{lem-de-fp-props}
  Consider any
  ML-FP (=DE-FP) $x=(x_w)\ne x^{**}$, for an anchored system, with left and right boundaries (if any) at $-\infty \le A < B \le \infty$.
  (Note that, if $A$ is finite, then necessarily, $x_w < 1$.)
  Then, $x$ satisfies the following additional properties. Function
    $x_w$  is strictly decreasing in $[A,B)$.
  The right derivative $x'^{+}_w$ exists at 
every point $w\in [A,B)$. The left derivative $x'^{-}_w$ exists, and  $x'^{-}_w \le x'^{+}_w$, at 
every point $w\in (A,B)$. (These derivatives are non-positive.)
Equation \eqn{eq-de-fp-new} holds for the right derivative $x'^{+}_w$ at every point $w\in [A,B)$.
The right derivative $x'^{+}_w$ is RCLL, the left derivative $x'^{-}_w$ is LCRL,
and each is negative, bounded away from $0$ on any finite interval $[A \wedge a,b]$, where $-\infty <a<b \le B$. 
  \end{lem}

  {\em Proof}  is essentially same as that of lemma 14 in \cite{St2021-coc}.

 \begin{lem}[Generalization of lemma 15 in \cite{St2021-coc}]
  \label{lem-hor-dist-no-local-max-new}
Consider any two non-degenerate (not equal to $x^{**}$) 
ML-FP (=DE-FP), $x$ and $\hat x$, for an anchored system, with positive speeds $v \le \hat v$ and arrival rates $\{\sigma_j\} \ge \{\hat \sigma_j\}$. 
(The systems may have different frame sizes, finite or infinite, in general.)  Denote by $\phi_u, ~u\in [0,1],$ and $\hat \phi_u, ~u\in [0,1],$ their corresponding inverse functions. 
Denote $b = \max\{x_\infty, \hat x_\infty\} < 1$, and consider the difference 
$\psi_u \doteq \hat \phi_u - \phi_u$ for $u\in (b,1]$.
Then, if for some  $\chi \in (b,1]$
$$
\psi_\chi = \max_{u\in [\chi,1]} \psi_u, 
$$
the difference $\psi_u$ is non-increasing in $(b,\chi]$. 
  \end{lem}

 Proof of Lemma~\ref{lem-hor-dist-no-local-max-new} is in Section~\ref{sec-monotone-diff}.

\begin{lem}[Extension of lemma 9 in \cite{St2021-coc}]
\label{lem-fsp-props-3}
 (i) Consider the anchored left-regulated system, with frame $[0,\infty)$. Consider
 the ML $x(\cdot)$, starting with the ``empty'' initial state $x(0)=x^\emptyset$. 
 Then $x(t)$ is monotone non-decreasing and $x(t) \to x^*$,
 where $x^*$ is an ML-FP, and moreover it is the {\em minimum ML-FP}. (That is, $x^* \le x$ for any ML-FP $x$.)
If $x^* \ne x^{**}$, then $x^*$ is proper. 

(ii) Consider the anchored right-regulated system, with frame $(-\infty,0]$. Consider
 the ML $x(\cdot)$, starting with the ``empty'' initial state $x(0)=x^\emptyset$. 
 Then $x(t)$ is monotone non-increasing and $x(t) \to x^*$,
 where $x^*$ is an ML-FP, and moreover it is the {\em maximum ML-FP}. (That is, $x^* \ge x$ for any ML-FP $x$.)
If $x^* \ne 1-x^{**}$, then $x^*$ is proper. (ML state $1-x^{**}$ corresponds to the entire particle mass at $-\infty$.)
  \end{lem}
  
Lemma~\ref{lem-fsp-props-3}(i) is lemma 9 in  \cite{St2021-coc}. The proof of Lemma~\ref{lem-fsp-props-3}(ii) is analogous.
We will use the common abbreviation ML-MFP for the minimum ML-FP in the left-regulated system and 
the maximum ML-FP in the right-regulated system -- it will be clear from the context which of the two we refer to.

\begin{lem}[Corollary of Lemma 20 in \cite{St2021-coc}]
\label{lem-bounded-first-moment-wave}
For an anchored two-side-regulated system, for any $v>0$,
there exists a unique ML-FP, 
and the dependence on $v$ is continuous. 
\end{lem}

\begin{lem}
\label{lem-de-fp-part}
Consider any unregulated on the left ML-FP (=DE-FP) $x$ in the anchored system. This ML may be regulated or not on the right; denote $\epsilon_*=0$ in the former case, and $\epsilon_* = x_{\infty} \ge 0$ in the latter. Denote $\delta = 1- x_{-\infty} \ge 0$. Fix any $\epsilon \in (\epsilon_*, 1-\delta)$, and let 
$w_*$ be the unique inverse of $x_w$ at $\epsilon$, i.e. $x_{w_*} = \epsilon$. Then, the renormalized trajectory
$$
\tilde x_{(-\infty,w_*]} = \frac{x_{(-\infty,w_*]} - \epsilon}{1-\epsilon-\delta}
$$
is a proper ML-FP ( = DE-FP) for the $(\epsilon,\delta)$-reduced system, right-regulated at boundary point $w_*$.
\end{lem}

Lemma~\ref{lem-de-fp-part} simply says that  if we take an ML-FP and truncate its domain to any left half-axis $(-\infty, w_*]$,
where $w_*$ is any point to the left of the right boundary (if any), we obtain (after renormalization) a proper ML-FP for a reduced system with 
the right boundary at $w_*$. The proof consists of the observation that if $x$ (being a DE-FP) satisfies \eqn{eq-de-fp-new} 
for the original system, then $\tilde x$ satisfies \eqn{eq-de-fp-new} for the reduced system; this observation, in turn, follows from the definition
of a reduced system.

By applying the notion of ML-MFP to a reduced system, we obtain, for example, the following property. Let the system parameters, and numbers $\epsilon \in [0,1)$, $w_* \in \R$, be fixed. If there exists a proper on the left ML-FP $x$ such that $x_{w_*} = \epsilon$, then there exists a proper on the left ML-FP, which dominates any other ML-FP in the interval $(-\infty,w_*]$. We will show later, in Theorem~\ref{thm-free-wave-unique-for-speed}, that a proper on the left ML-FP $x$, such that $x_{w_*} = \epsilon$,
is in fact unique.

\section{Preliminary definitions and results (generalizing those in \cite{St2021-coc}), related to SSAI}
\label{sec-res-coc-gen}

This section contains definitions and results, generalizing those in \cite{St2021-coc}, which are related to SSAI.
When necessary, we provide proofs, deferred to Section~\ref{sec-proofs-new}.

For a given $n$ consider the Markov process $x^n(\cdot)$ for an anchored system,
and the corresponding centered process $\mathring x^n(\cdot)$.
Process $\mathring x^n(\cdot)$ is {\em not} Markov, because the information about the location of the boundaries is ``lost'' after centering;
it is a projection of Markov process $x^n(\cdot)$.

\begin{thm}[Theorem 3 in \cite{St2021-coc}]
\label{thm-bounded-first-moment}
Consider the anchored left-regulated system with frame $[0,\infty)$.
Suppose Assumptions~\ref{cond-second-moment},  \ref{cond-ntriv-add-ii} and \ref{cond-ntriv-add-i}  hold.
Then there exist $\bar C>0$ and $\bar n$ such that, uniformly in $n \ge \bar n$ and corresponding (for each $n$) speeds $v$ %$\lambda^n$
such that the process $x^n(\cdot)$ is stable, we have
\beql{eq-tight-uniform}
\E \Phi_1(\mathring x^n(\infty)) \le \bar C.
\eeql
\end{thm}

Proof of Theorem~\ref{thm-bounded-first-moment} is in Section~\ref{sec-closeness-proof}. It is an extended and clarified version of the proof in \cite{St2021-coc}. It is provided for completeness, and because versions of the arguments in this proof are used in
other parts of this paper.
We note that the result analogous to Theorem~\ref{thm-bounded-first-moment} can be obtained for the right-regulated system with frame $(-\infty,0]$
as well, with analogous proof.

Note that for an unregulated (free) system the centered process $\mathring x^n(\cdot)$ 
{\em is} Markov. 
This process only changes its state (jumps) upon job arrivals. Therefore, 
the speed parameter $v$ is irrelevant, and 
scaling of the job arrival rate 
$\lambda^n$ only scales the rate at which changes take place. Consequently, changing $v$ and/or $\lambda^n$ (as long as it is positive)
  has no impact on the existence/non-existence/form of centered process $\mathring x^n(\cdot)$ stationary distribution in free system.

A version of Theorem~\ref{thm-bounded-first-moment} holds for the free system, 
under an additional technical condition needed to be able to demonstrate
stability (positive Harris recurrence).
State $\mathring x^n(t)$ of the free system can be equivalently described as $p^n(t) = (w_1(t), \ldots, w_n(t)) \in \R^n$,
where $w_1, w_2, \ldots, w_n$ are the locations of the $n$ particles (not necessarily ordered in any way)
w.r.t. the mean $\bar x^n(t)$.
A subset $B \subset \R^n$ is called {\em small} (see, e.g. \cite{Bramson-book}) for the process $p^n(\cdot)$,
 if there exists $\tau >0$ and a finite measure 
$\mathcal M$ on $\R^n$, such that, uniformly on the initial state $p^n(0)\in B$, the distribution of $p^n(\tau)$ dominates $\mathcal M$.
(Existence of a small set allows one to use Nummelin splitting to view the process as having a renewal atom.) Note that,
for process $p^n(\cdot)$,
 the property of a set being small
does not depend on $\lambda^n$, as long as $\lambda^n>0$.

\begin{assumption}
\label{assump-small}
In the free system, for all large $n$, 
\beql{eq-small-set2}
\mbox{for any $C>0$ the (closed) set $\{|p^n| \le C\}$ is small for $p^n(\cdot)$.} 
\eeql
\end{assumption}

Assumption~\ref{assump-small} holds in many cases of interest. For example, it is clearly satisfied when process $p^n(\cdot)$
happens to be an irreducible countable Markov chain and sets $\{|p^n| \le C\}$ are finite. 
As another example, Assumption~\ref{assump-small}
 is satisfied in the case when a slightly stronger version of Assumption~\ref{cond-ntriv-add-ii} holds, as shown in the following

\begin{lem}
\label{lem-small-sufficient}
Suppose there exists a class $j$, $0\le m_j < k_j$, and $\epsilon>0$, such that the component size distribution $F_j$ dominates a scaled 
(by a positive constant) measure such that $m_j$ components are zero and the remaining $d_j - m_j$ components distribution is the
Lebesgue measure on $[0,\epsilon]^{d_j - m_j}$. Then Assumption~\ref{assump-small} holds.
\end{lem}

Proof of Lemma~\ref{lem-small-sufficient} is in Section~\ref{sec-proof-small-sufficient}.

\begin{thm}[Theorem 11 in \cite{St2021-coc}]
\label{th-closeness-free}
Consider a free system. 
Suppose Assumptions~\ref{cond-ntriv-add-ii},  \ref{cond-ntriv-add-i}, and \ref{assump-small} hold.
Then, there exists  $\bar n$ such that for all $n \ge \bar n$,  $\mathring x^n(\cdot)$ 
is positive recurrent.%\footnote{In fact, Assumption~\ref{cond-ntriv-add-i} is {\em not} needed for positive recurrence. See Appendix A in https://arxiv.org/abs/2507.11779}
\footnote{In fact, Assumption~ \ref{cond-ntriv-add-i} is {\em not} needed for positive recurrence. See Theorem~\ref{th-free-stabil-gen} in Appendix~\ref{app1}.}
If, in addition, Assumption~\ref{cond-second-moment} holds, then there exist $\bar C>0$ such that
\beql{eq-tight-uniform-free}
\E \Phi_1(\mathring x^n(\infty)) \le \bar C.
\eeql
\end{thm}

Proof of Theorem~\ref{th-closeness-free} is an adjustment of the proof of Theorem~\ref{thm-bounded-first-moment}; the necessary adjustments are specified in \cite{St2021-coc}. We note that the positive recurrence proof in \cite{St2021-coc} does {\em not} use Assumption~\ref{cond-second-moment}.

If the process $\mathring x^n(\cdot)$ in a free system is positive recurrent, then the average velocity $v_n$ of the unanchored free particle process 
$x^n(\cdot)$ is well-defined, with multiple possible equivalent definitions. For example, 
\beql{eq-vn-def-wp1}
v_n = \lim_{t\to\infty} \Delta^n(t;\nu)/t, ~~w.p.1,
\eeql
for any quantile $\nu\in (0,1)$, and any fixed initial state $x^n(0)$, where $\Delta^n(t;\nu)$ is the random displacement of the $\nu$-quantile of $x^n(t)$ in the time interval $[0,t]$. Another definition is:
\beql{eq-vn-def-stat}
v_n = \E \Delta^n(t;\nu)/t,
\eeql
for any quantile $\nu\in (0,1)$ and time $t>0$, 
when $\mathring x^n(\cdot)$ is in the stationary regime.

\begin{condition}[DE-FP/ML-FP uniqueness]
\label{cond-unique}
Consider an anchored left-regulated system, with frame $[0,\infty)$.
We say that, in the anchored left-regulated system with frame $[0,\infty)$, 
this condition holds for speed $v > 0$, if 
for any sufficiently small $b \ge 0$, 
there is at most one DE-FP $x$ with $x_\infty = b$. 
\end{condition}

\begin{prop}[Key part of the proof of theorem 4 in \cite{St2021-coc}]
\label{prop-ssai-key}
Consider the left-regulated system with the frame $[0,\infty)$.
Suppose Assumptions~\ref{cond-second-moment}, \ref{cond-ntriv-add-ii} and \ref{cond-ntriv-add-i} hold.
Consider our asymptotic regime, with the number of particles $n\to\infty$,
such that for all large $n$ the process is stable and, moreover, $\rho^n = \E x^n_0(\infty) \to \rho < 1$.
Suppose that Condition~\ref{cond-unique} holds for speed $v>0$.
Then SSAI holds for this $v$, i.e. $x^n(\infty) \Rightarrow x^{*}$, where $x^{*}$ is the ML-MFP and it is such that $x^{*}_0=\rho$.
\end{prop}

Consider a system unregulated on the right.
Consider a fixed job class $j$.  
Let $(Z_1, \ldots, Z_{d_j})$ be the ordered vector of particle locations in
 a selection set (of a class $j$ job); namely, $Z_i$ are the locations of selected particles,
ordered to form a non-decreasing sequence, $Z_{i} \le Z_{i+1}$. Denote by $D=(D_1, \ldots, D_{d_j})$
the location-differential vector, where $D_1 =0$, $D_i = Z_{i} - Z_{i-1}, ~i=2,\ldots,d_j$. 
Denote by $\eta_j(D)$ the (random) total displacement of all selected particles due to
a class $j$ job arrival, given the selection set location-differential vector is $D$.
It is easy to see that 
\beql{eq-eta-cont}
\mbox{$\E \eta_j(D)$ is continuous in $D$.}
\eeql

\begin{definition}[$D$-monotonicity]
\label{def-d-monotone}
A class $j$  
is called $D$-monotone decreasing, if its distribution $F_j$ of component sizes is such that $\E \eta_j(D)$ is non-increasing in $D$, i.e.
$D \ge D'$ implies $\E \eta_j(D) \le \E \eta_j(D')$.
\end{definition}

\begin{lem}
\label{lem-D-monotone-subclasses}
If a class $j$ is $D$-monotone decreasing, then so are all its subclasses $j^{m,\ell}$.
\end{lem}

The proof of Lemma~\ref{lem-D-monotone-subclasses} is straightforward. Consider any finite locations of $m-\ell$ particles in the reduced system;
consider original system with $m-\ell$ particles located as in the reduced system, $\ell$ particles located at $-C$, $d_j-m$ particles located at $C$;
let $C \uparrow \infty$.

\begin{definition}[Increasing hazard rate (IHR)]
\label{def-ihr}
A distribution $H(\cdot)$ on $\R_+$ has increasing hazard rate (IHR), if
\beql{eq-ihr-cond}
\frac{H(y_1+\Delta) - H(y_1)}{1-H(y_1)} \le \frac{H(y_2+\Delta) - H(y_2)}{1-H(y_2)}, ~~\forall ~ 0 \le y_1 \le y_2 ~\mbox{such that $H(y_2)<1$}, ~\forall \Delta \ge 0.
\eeql
\end{definition}

Examples of IHR distributions include: the exponential distribution (which, in fact, has constant hazard rate); deterministic distribution 
concentrated on a single point $a \ge 0$; uniform distribution on a finite interval $[0,a]$. Also, if a random variable $A$ has an IHR distribution
and $a\ge 0$ is a constant, then the truncated random variable $A \wedge a$ also has an IHR distribution.

We will say that a distribution $F_j$ is a mixture of distributions $F_{j,\chi}$, parameterized by $\chi$,
if it is obtained by averaging $F_{j,\chi}$ w.r.t. some probability distribution on the values of $\chi$.
(In other words, a realization of $F_j$ is obtained by first choosing a realization of $\chi$ according its distribution
and then choosing a realization of $F_{j,\chi}$.)

\begin{prop}[Theorem 6 in \cite{St2021-coc}]
\label{thm-main-D-monotonicity}
(i) If all classes $j$ 
are D-monotone decreasing,
  then Condition~\ref{cond-unique} holds for any speed $v > 0$.
  
(ii) Suppose the distribution $F_j$ for a class $j$ is such that the component sizes are i.i.d. with an IHR distribution $H_j(\cdot)$. 
(Or, more generally, $F_j$ may be a mixture of such distributions.)
Then class $j$ is 
$D$-monotone decreasing.
\end{prop}

\section{Main results}
\label{sec-main-results-new}

\begin{thm}
\label{thm-bounded-first-moment-tws}
Consider ML-FPs for a family of anchored systems, with the following set of parameters taking values within a compact set: speed $v$; 
limiting arrival rates $\{\sigma_j\}$; and the distributions $\{F_j\}$ of the job component sizes
(convergence of distributions is in the weak sense). 
Assume further that the family of parameters is such that: speed $v$ is bounded away from $0$; Assumption~\ref{cond-second-moment} holds uniformly;
Assumption~\ref{cond-ntriv-add-ii} holds uniformly in that it holds for the same class $j$, with the corresponding probability bounded away from $0$,
and $\sigma_j$ bounded away from $0$. Then, for some $\bar C>0$, uniformly in that family of parameters and
uniformly across all free and regulated systems, an ML-FP $x$ is such that
$\Phi_1(\mathring x)\le \bar C < \infty.$ 
\end{thm}

Proof of Theorem~\ref{thm-bounded-first-moment-tws} is in Section~\ref{sec-proof-bounded-first-moment}.

\begin{thm}
\label{thm-free-wave-exists}
Suppose Assumptions~\ref{cond-second-moment} and \ref{cond-ntriv-add-ii} hold. 
Then, for the anchored free system, a proper ML-FP with a positive speed exists.
\end{thm}

{\em Proof.} For a two-side-regulated system, with frame $[0,B]$ there exists a unique ML-FP   
for any speed $v>0$. For a sequence of frame sizes $B \to \infty$, choose drift $v(B)$, such that the median of the fixed point is exactly in the middle, at point $B/2$. Let $x^{(B)}$ be such ML-FP shifted left by $B/2$, so that the median is at $0$. Let us choose a subsequence of $B$ along which $v(B) \to v$. Limiting speed $v$ cannot be $0$ or $+\infty$, because otherwise the median could not stay in the middle. (Indeed, if the limiting speed is infinity, then the probability mass of the left boundary has to stay bounded away from $0$ -- otherwise
the average drift to the right cannot match the average drift to the left; in fact, that mass goes to $1$, so the distribution concentrates on the left boundary.
If the limiting speed is zero, then the distribution, in a weak limit, must concentrate at the right boundary  -- otherwise
the average drift to the right cannot go to $0$ to match the average drift to the left.) Then, by Theorem~\ref{thm-bounded-first-moment-tws}, 
$\Phi_1(\mathring x^{(B)}) \le C <\infty$, uniformly in all $c$. This in turn implies that $\Phi_1(x^{(B)}) \le C <\infty$, uniformly in all $B$,
because $|\bar x^{(B)}|$ must stay bounded.
So, we can choose a further subsequence of $B$ along which $x^{(B)} \to x$,
where $x$ is proper. Then, by Lemma~\ref{lem-cont}, $x$ is a proper free ML-FP with speed $v$. 
$\Box$

The next Theorem~\ref{thm-free-wave-unique-for-speed} is one the key results of the paper. 
It shows that for a given speed $v>0$ the proper left tail of a DE-FP (for an anchored left-unregulated  system) is determined uniquely.
This leads, for a given speed $v>0$, to the uniqueness of a ML-FP  
with proper left tail,
for any system unregulated on the left, which may or may not be regulated on the right. To state Theorem~\ref{thm-free-wave-unique-for-speed} we need some additional definitions.

  Consider an artificial $(\epsilon,0)$-reduced system, with $\epsilon\ge 0$, and with frame $(-\infty,0]$, 
where the arrival rate per server is $\alpha = \lambda \sum_j \pi_j d_j = \sum_j \sigma_j$
and each job having a single component with distribution 
$$
H(w) = \sum_j \frac{\pi_j d_j}{\sum_\ell \pi_\ell d_\ell} H_j(w),
$$
where $H_j(w)$ is the marginal CDF of $\xi_1^{(j)}$. Consider the following exponential function 
\beql{eq-chi-def}
\chi^{(\epsilon)}_w = 1- (1-\epsilon) e^{\beta w}, ~w\le 0
\eeql
There is a unique parameter $\beta>0$, such that this function is a DE-FP for the artificial system with speed $v$, 
as long as 
\beql{eq-cond-left-tail-exists}
\alpha \int_0^{\infty} w dH(w) > v.
\eeql
Indeed, 
let $L(\beta) = \int_0^{\infty} e^{-\beta z} d H(z), ~\beta\ge 0,$ be the Laplace transform of distribution $H(\cdot)$, and 
$\bar L(\beta) = \int_0^{\infty} e^{-\beta z} \bar H(z) dz = [1-L(\beta)]/\beta$, $\bar H(w) = 1- H(w)$. Then 
$\chi^{(\epsilon)}$ is a DE-FP if
$$
\int_{-\infty}^z [-\chi^{(\epsilon)}_w]' \alpha \bar H(z-w) dw = -[\chi^{(\epsilon)}_z]' v, ~~z \le 0,
$$
which reduces to
\beql{eq-cond-beta}
\alpha \bar L(\beta) = v.
\eeql
The LHS is convex, decreasing to $0$ as $\beta\to\infty$, with $\bar L(0) =  \int_0^{\infty} w dH(w)$, which is the mean component size.
So, \eqn{eq-cond-beta} uniquely determines $\beta$ as long as \eqn{eq-cond-left-tail-exists} holds.

Note that for any $0 \le \epsilon < 1$, 
$$
(\chi^{(\epsilon)}_w-\epsilon)/(1-\epsilon) = \chi^{(0)}_w= 1-e^{\beta w}, ~w\le 0,
$$
and $\chi^{(0)}$ is nothing else, but the stationary distribution of a particle moving in $(-\infty,0]$,
with the rate $\alpha$ Poisson process of independent jumps to the right (but not over the right boundary at $0$), 
with size distribution $H(\cdot)$, and moving left at speed $v$ between the jumps.
We conclude that $\chi^{(\epsilon)}$ is the unique ML-FP (=ML-MFP) of the artificial $(\epsilon,0)$-reduced system
with the right boundary at $0$.

As will be seen in the proof of Theorem~\ref{thm-free-wave-unique-for-speed}, $\chi^{(\epsilon)}$ gives the asymptotic shape of the left tail of any DE-FP with speed $v$.
Informally, this is because, when $1-\epsilon$ is close to $0$, the $(\epsilon,0)$-reduced original system is approximately the $(\epsilon,0)$-reduced artificial system.

\begin{thm}
\label{thm-free-wave-unique-for-speed}
(i) In the anchored system, for any $v$ satisfying \eqn{eq-cond-left-tail-exists}, there exists the unique DE-FP with the proper left tail.
This solution will either (a) hit the $w$-axis or (b) will converge to some $\epsilon_* \ge 0$ as $w\uparrow \infty$.
In case (a): for any $0 \le \epsilon <1$ there is a unique ML-FP with proper left tail for the $(\epsilon,0)$-reduced right-regulated system.
In case (b): for any $\epsilon_* < \epsilon <1$ there is a unique ML-FP with proper left tail for the $(\epsilon,0)$-reduced right-regulated system;
for $\epsilon = \epsilon_*$ there is a unique ML-FP with proper left and right tails for the $(\epsilon,0)$-reduced free system;
for $\epsilon < \epsilon_*$ there is a unique ML-FP with proper left tail, but not proper on the right, for the $(\epsilon,0)$-reduced free system.

(ii) In the anchored system, for any $v$ satisfying \eqn{eq-cond-left-tail-exists}, and any $0 \le \epsilon <1$,
there exists at most one ML-FP with proper left tail for the $(\epsilon,0)$-reduced system
with speed $v$; this ML-FP is either for the right-regulated system or for the free system, with proper or not proper right tail.
\end{thm}

Proof of Theorem~\ref{thm-free-wave-unique-for-speed} is in Section~\ref{sec-free-wave-unique-for-speed}.

\begin{thm}
\label{thm-free-wave-spectrum} 
Suppose Assumptions~\ref{cond-second-moment} and \ref{cond-ntriv-add-ii} hold.

(i) For the anchored free system, there exists a range $[v_{min},v_{max}]$ of speeds $v$, with $0 < v_{min} \le v_{max} < \infty$ for which a proper ML-FP $x$ exists. 

(ii) If $x$ and $\hat x$ are two ML-FPs with speeds $v < \hat v$, respectively, and $\phi$ and $\hat \phi$ their respective inverses, 
then the difference $\hat \phi_u - \phi_u, ~ 0 \le u \le 1$, is non-increasing.

(iii) For any $v < v_{min}$,  
the unique ML-FP with proper left tail is an ML-FP for a right-regulated system. 
For $v > v_{max}$, the unique ML-FP with proper left tail is an ML-FP for the free system, with non-proper right tail.
\end{thm}

{\em Proof.} (i) At least one proper ML-FP exists by Theorem~\ref{thm-free-wave-exists}. Let $v_{min} \ge 0$ be the infimum of speeds, for which a proper ML-FP exists. 
By Assumption~ \ref{cond-finite-mean}(ii), upon an arrival of a job, uniformly on 
the mutual positions of the particles it selects, the expected total (positive) displacement of the selected particles 
is bounded away from $0$ by a constant $\delta>0$. Then, a proper ML-FP speed must be at least $v\ge \lambda \delta$.
Therefore, $v_{min} > 0$. Also, upon an arrival of a job, uniformly on 
the mutual positions of the particles it selects, the expected total (positive) displacement of the selected particles 
is upper bounded by the expected total size of all components, which is uniformly bounded by a finite $\delta'$.
Therefore, $v_{max} \le \lambda \delta' < \infty$. There must be ML-FPs for both speeds $v_{min}$ and $v_{max}$,
by Lemma~\ref{lem-cont} and Theorem~\ref{thm-bounded-first-moment-tws}. 

(ii) To see the monotonicity of the difference of inverses, pick any $\epsilon \in (0,1)$. Let as shift $x$ and $\hat x$ so that $x_0=\hat x_0 =\epsilon$.
In the interval $(-\infty,0]$, $x$ and $\hat x$ are the scaled versions of ML-MFPs for the $(\epsilon,0)$-reduced system with frame $(-\infty,0]$,
 for the speeds
$v$ and $\hat v$, respectively. Therefore, in $(-\infty,0]$, $x$ dominates $\hat x$. Then, by Lemma~\ref{lem-hor-dist-no-local-max-new},
the difference $\hat \phi_u - \phi_u, ~ 0 \le u \le 1$, is non-increasing for $0 \le u \le \epsilon$. Since $\epsilon$ can be arbitrarily close to $1$,
the difference is non-increasing for $0 \le u \le 1$.

Now, if ML-FPs exist for speeds $v$ and $\hat v$, with $v < \hat v$, then the ML-FP exists
for any speed $v' \in (v, \hat v)$. Indeed, if we consider the point where the former two intersect,
then there exists a unique DE-FP solution with speed $v'$ with proper left tail, which contains that intersection point;
both to the left and to the right of the intersection point, this solution is sandwiched between the former two ML-FPs,
and therefore has to be proper on the right as well.

(iii) If $v < v_{min}$, then \eqn{eq-cond-left-tail-exists} is satisfied, and
the ML-FP with proper left tail is the unique DE-FP solution $x$, proper on the left. Let us shift $x$ solution so that
it intersects with ML-FP $\hat x$ for speed $v_{min}$ at a given point, say at $w=0$. Then, by (ii), the $x_w$ is upper bounded by $\hat x_w$ 
for $w\ge 0$. Then, $x$ must ``hit'' $w$-axis (i.e., we must have case (a) in Theorem~\ref{thm-free-wave-unique-for-speed}(i)), because otherwise 
$x$ would be a proper free ML-FP with speed less than $v_{min}$, which contradicts the $v_{min}$ definition. Therefore, $x$ is 
the ML-FP for the system with frame $(-\infty,w_*]$, where $w_*$ is point, where $x$ hits
$w$-axis. Similarly, for $v > v_{max}$, the unique DE-FP solution $x$ must be such that $x_\infty = \epsilon_* >0$
(i.e., case (b) in Theorem~\ref{thm-free-wave-unique-for-speed}(i)); otherwise, comparing $x$
to the proper ML-FP for the speed $v_{max}$, we would obtain a contradiction to the definition of $v_{max}$.
$\Box$

Consider an ML $x(\cdot)$ in unanchored system, with proper initial state $x(0) \in \cx^{pr}$. 
We say that the velocity of $x(\cdot)$ is lower bounded by $v \ge 0$ if
$$
\liminf_{t\to\infty} \Delta(t;\nu)/t \ge v,
$$
for any $\nu\in (0,1)$,  where $\Delta(t;\nu)$ is the displacement of the $\nu$-quantile of $x(t)$ in the time interval $[0,t]$.

\begin{lem}
\label{lem-right-regulated-wave-as-lower-bound111} 
Suppose there exists a right-regulated ML-FP with speed $v$ in the anchored system with arrival rates $\{\tilde \sigma_j\}$ dominated by 
the rates $\{\sigma_j\}$ in the original system.
Then the velocity of the ML $x(\cdot)$ in the unanchored free system, starting from initial state $x(0) = x^\emptyset$, is lower bounded by $v$.
\end{lem}

The proof of Lemma~\ref{lem-right-regulated-wave-as-lower-bound111} follows by coupling and monotonicity, if WLOG we assume that 
the initial location of the right boundary is at $0$.

\begin{thm}
\label{thm-right-regulated-wave-as-lower-bound}
Suppose for any speed $\tilde v < v$, 
there exists a right-regulated ML-FP with speed $\tilde v \le v$ in an anchored system with arrival rates $\{\tilde \sigma_j\}$,
such that $\{\tilde \sigma_j\} \le \gamma \{\sigma_j\}$, where $\gamma<1$ and $\{\sigma_j\}$ are the rates in the original system.
Then the velocity of any proper ML $x(\cdot)$ in the unanchored free system, starting from any initial state $x(0) \in \cx^{pr}$, is lower bounded by $v$.
\end{thm}

{\em Proof.} 
Let us fix a small $\delta>0$, and assume WLOG that the initial state $x(0)$ of ML $x(\cdot)$ is such that $x_0(0)=1-\delta$. 
Then $x(\cdot)$ dominates ML $x^\ell(\cdot)$ with initial state $x^\ell(0) = (1-\delta) x^\emptyset \le x(0)$; this initial state is such that mass $\delta$ is located at $-\infty$
and mass $1-\delta$ is located at $0$. We see that the velocity of $x^\ell(\cdot)$, if we only look at $\nu$-quantiles with $\nu \in (\delta,1)$, 
is same as the velocity of the ML starting from $x^\emptyset$ in the $(0,\delta)$-reduced system. 
Now note that in the $(0,\delta)$-reduced system the job class arrival rates $\{\hat \sigma_j\}$ are only slightly different from $\{\sigma_j\}$,
and therefore $\{\tilde \sigma_j\} \le  \{\hat \sigma_j\}$.
Applying Lemma~\ref{lem-right-regulated-wave-as-lower-bound111}, we see that the velocity of $x^\ell(\cdot)$,
limited to $\nu$-quantiles with $\nu \in (\delta,1)$,
 is at least some $\tilde v <v$, which can be made arbitrarily close to $v$ by making $\delta$ small. Recalling that $x(\cdot)$ dominates $x^\ell(\cdot)$,
 we conclude that the velocity of $x(\cdot)$ is lower bounded by $v$. 
$\Box$

For some of the results we will need the following
\begin{assumption}
\label{assump-101}
All classes $j$ and their sub-classes of the form $j^{m,0}$ (i.e., those which are present in an $(\epsilon,0)$-reduced system with $\epsilon>0$), have strictly positive arrival rates $\sigma_j$ in the original model. 
(Note that subclasses $j^{m,0}$ are exactly the subclasses considered in \cite{St2021-coc}.)
\end{assumption}

\begin{thm}
\label{thm-free-wave-speed-as-lower-bound}
Suppose Assumption~\ref{assump-101} holds.
Suppose $x^*$ is a proper ML-FP with speed $v$ in the anchored free system. 
Then the velocity of any ML $x$ in the unanchored free system, with proper initial state $x(0) \in \cx^{pr}$, is lower bounded by $v$.
\end{thm}

{\em Proof.}  
Let us fix a small $\epsilon>0$. Consider the ML-FP $x^*$ such that $x^*_0=\epsilon$. Then, $x^*$ in the interval $w \le 0$ is the unique ML-FP 
for the $(\epsilon,0)$-reduced system with the right boundary at $0$. Since $\epsilon$ is small, the job class arrival rates $\{\hat \sigma_j\}$ 
in the reduced system
are close to the job class arrival rates $\{\sigma_j \}$  
in the original system. Since all $\sigma_j $ are positive
(by Assumption~\ref{assump-101}), 
in the reduced system
we can slightly reduce the speed from $v$ to $\tilde v$, and reduce $\lambda$ by the same factor
to $\tilde \lambda = \lambda (\tilde v /v)$, so that the job class arrival rates $\{\tilde \sigma_j = \hat \sigma_j (\tilde v /v)\}$  
in the reduced system are such that $\{\tilde \sigma_j\} < \{\sigma_j \}$.
Note that $\tilde v$ can be arbitrarily close to $v$ if 
$\epsilon$ is small enough. Application of Theorem~\ref{thm-right-regulated-wave-as-lower-bound} completes the proof.
$\Box$

\begin{thm}
\label{thm-free-wave-unique-cond}
For the anchored free system, each of the following conditions (i) and (ii) is sufficient for the uniqueness of an ML-FP speed (and then for the existence and uniqueness of a ML-FP):\\
 (i) Assumption~\ref{assump-101}; \\
 (ii) Assumptions~\ref{cond-second-moment} and \ref{cond-ntriv-add-ii}, and D-monotonicity.
\end{thm}

{\em Proof of Theorem~\ref{thm-free-wave-unique-cond}.} 
(i) By Theorem~\ref{thm-free-wave-speed-as-lower-bound}, the velocity of any unanchored free ML is lower bounded by the speed $v$ of any ML-FP. Therefore, the existence of two free traveling waves with different speeds is impossible.

(ii) Suppose $x$ and $\hat x$ are two different ML-FP, with speeds $v < \hat v$. 
By Theorem~\ref{thm-free-wave-spectrum} the difference of the inverses $\hat \phi_u - \phi_u$ is non-increasing.
Consider relation \eqn{eq-speed-equal-workload-arrivals} for $x$ and $\hat x$: the LHS has to be strictly greater for $\hat x$; on the other hand,
the RHS has to be not greater for $\hat x$ by D-monotonicity (and $\hat \phi_u - \phi_u$ being non-increasing).
The contradiction shows that $v < \hat v$ is impossible.
$\Box$

\begin{thm}
\label{thm-free-wave-vmin-vmax-monotone}
Consider the free system. 

(i) Consider two sets of parameters,  with arrival rate vectors $\{\sigma_j\} \le \{\hat \sigma_j\}$ and with component size distributions $F_j$ in the former set being stochastically dominated by corresponding distributions $\hat F_j$ in the latter. Then, we have
$v_{min} \le \hat v_{min}$ and $v_{max} \le \hat v_{max}$, where $\hat v_{min}$ and $\hat v_{max}$ are the minimum and maximum ML-FP speeds 
in the latter system.

(ii) Suppose, in addition, that $\hat \sigma \downarrow \sigma$, and component size distributions converge, and Assumptions~\ref{cond-second-moment} and \ref{cond-ntriv-add-ii} hold uniformly for all systems involved. Then
$\hat v_{max} \downarrow v_{max}$ and the corresponding ML-FPs for those speeds converge as well, $\hat x^{max} \to x^{max}$. Analogously, in the case when $\tilde \sigma \uparrow \sigma$, and other parameters converge from below, we have $\tilde v_{min} \uparrow v_{min}$ and 
$\tilde x^{min} \to x^{min}$.
\end{thm}

{\em Proof of Theorem~\ref{thm-free-wave-vmin-vmax-monotone}.} 
(i) Fix any $v < v_{min}$. Consider the unique proper on the left, right-regulated ML-FP $x$ for $\{\sigma_j\}$. This is a DE-FP, hitting $w$-axis at some point. But then, for the same speed $v$, the proper on the left DE-FP $\hat x$ for $\{\hat \sigma_j\}$ , intersecting $x$ at a specific point is (by monotonicity 
of inverses -- see Lemma~\ref{lem-hor-dist-no-local-max-new}) is such that it is below $x$ to the right of the intersection point. This means that $\hat x$ also hits $w$-axis. Then $\hat v_{min} \ge v$, and then $\hat v_{min} \ge  v_{min}$.

Consider speed $v_{max}$ and the corresponding unique proper ML-FP $x$ for $\{\sigma_j\}$. The proper on the left DE-FP $\hat x$ for 
$\{\hat \sigma_j\}$, intersecting $x$ at a specific point is (by monotonicity of inverses -- see Lemma~\ref{lem-hor-dist-no-local-max-new}) is such that it is below $x$ to the right of the intersection point. So, it cannot be non-proper on the right. Then $\hat v_{max} \ge v_{max}$.

(ii) The second statement follows by Lemma~\ref{lem-cont}, along with the uniform boundedness of the first moment (by Theorem~\ref{thm-bounded-first-moment-tws}).
$\Box$

\begin{thm}
\label{thm-free-wave-to-regulated-unique} 
Suppose Assumptions~\ref{cond-second-moment} and \ref{cond-ntriv-add-ii} hold. 

(i) If free ML-FP speed range is $[v_{min},v_{max}]$, then for the anchored left-regulated system with frame $[0,\infty)$,
for any $v>v_{max}$, there exists a unique proper ML-FP (which is then equal to the ML-MFP), and Condition~\ref{cond-unique} holds.

(ii) As $v \downarrow v_{max}$, the corresponding proper left-regulated ML-FP $x^*$ is such that its load 
$x^*_0 \uparrow \bar \rho \le 1$.

(iii) If $\bar \rho <1$, then the following additional properties hold.
We have convergence $x^* \uparrow x^{*,max}$, where $x^{*,max}$ is the ML-MFP for speed $v_{max}$, and $x^{*,max}_0=\bar \rho$.
Furthermore, for the speed $v_{max}$ there is the entire family of left-regulated proper ML-FP $\tilde x^*$ with  $\tilde x^*_0 \in [\bar \rho, 1)$,
and with $\tilde x^* \to x^{*,max,free}$ as $\tilde x^*_0 \uparrow 1$, where $x^{*,max,free}$ is the unique free ML-FP for speed $v_{max}$.
(To be precise, the convergence $\tilde x^* \to x^{*,max,free}$ holds if we take versions of $\tilde x^*$ shifted to intersect $x^{*,max,free}$
at a given point -- say, $\tilde x^*_0 \to x^{*,max,free}_0$.)
Any two ML-FL $\tilde x^*$ and $\hat x^*$ from the above family, with $\tilde x^*_0 < \hat x^*_0$, are such that the difference of their 
inverses $\hat \phi_u - \tilde \phi_u$ is  non-increasing.

(iv) If $\bar \rho = 1$, then we have convergence of (appropriately shifted versions of) $x^*$ to the unique free ML-FP 
$x^{*,max,free}$ for speed $v_{max}$.

(v) If Condition~\ref{cond-unique} holds for all $v$, then $\bar \rho = 1$.
\end{thm}

{\em Proof of Theorem~\ref{thm-free-wave-to-regulated-unique}.} 
(i) Consider any speed $v>v_{max}$. Then, the ML-FP $x$, unregulated on the left, with proper left tail, is unique and not proper, having $x_\infty = \epsilon' >0$. For this fixed speed $v$, consider a family of left-regulated ML-FP $y$, which all intersect with $x$ at $w=0$, with the left boundary moving to the left to $-\infty$. 
The limit of $y$ must be exactly $x$. Using this fact, we see that at some point will obtain a left-regulated ML-FP $y$, which is not proper on the right -- otherwise, $\Phi_1(\mathring y)$ would go to infinity, which is impossible. So, there exists a left-regulated ML-FP with speed $v$, not proper on the right, i.e., $y_\infty =\epsilon>0$. Then, by the same argument as in the proof of theorem 8 in \cite{St2021-coc}, we conclude that there exists a continuous family of DE-FP solutions in frame $[0,\infty)$, with speed $v$, and with asymptotic values on the right continuously changing from $0$ to $\epsilon$. This proves that the left-regulated proper ML-FP is unique (and then equal to the ML-MFP), and that Condition~\ref{cond-unique} holds.

(ii) Since $x^*$ for $v>v_{max}$ is the ML-MFP, $x^*_0$ is non-increasing with $v$, by monotonicity.

(iii) The existence of the limiting ML-FP $x^{*,max}$ with specified properties easily follows from Lemma~\ref{lem-cont}.
Consider $x^{*,max}$ with the left boundary at $0$, and $x^{*,max,free}$ intersecting $x^{*,max}$ at $w=0$ (i.e., $x^{*,max,free}_0=x^{*,max}_0$).
The specified family $\tilde x^*$ of left-regulated ML-FP is obtained by considering the family of DE-FP with the left boundary moving left to $-\infty$,
and intersecting $x^{*,max,free}$ and $x^{*,max}$ at the same point at $w=0$. 
The difference of inverses $\hat \phi_u - \tilde \phi_u$ is  non-increasing for any two ML-FP in this family,
including $x^{*,max,free}$ and $x^{*,max}$ as its extreme points, is by Theorem~\ref{thm-free-wave-spectrum}(ii).
Then any $\tilde x^*$ is sandwiched between $x^{*,max,free}$ and $x^{*,max}$
both in $(-\infty,0]$ and in $[0,\infty)$, and therefore each $\tilde x^*$ is proper. 

(iv) This easily follows if we consider the shifted versions of $x^*$, intersecting $x^{*,max,free}$ at the same point at $w=0$.

(v) If $\bar \rho < 1$ would hold, then by (iii) we would have multiple left-regulated proper ML-FP for the same speed $v_{max}$.
This would contradict Condition~\ref{cond-unique}.
$\Box$

From Theorem~\ref{thm-free-wave-to-regulated-unique} and Proposition~\ref{prop-ssai-key} we obtain the following result on the SSAI
for left-regulated systems,
Theorem~\ref{thm-free-wave-to-ssai}, 
which generalizes and unifies main SSAI results 
 in \cite{St2021-coc}, as will be seen form its Corollaries~\ref{thm-free-wave-to-ssai-cor1} and \ref{thm-free-wave-to-ssai-cor2}.

\begin{thm}
\label{thm-free-wave-to-ssai}
Suppose Assumptions~\ref{cond-second-moment}, \ref{cond-ntriv-add-ii} and \ref{cond-ntriv-add-i} hold.
If free ML-FP speed range is $[v_{min},v_{max}]$, then SSAI 
holds for the left-regulated systems with frame $[0,\infty)$ for speeds $v>v_{max}$, 
which correspond to ML-FP loads within $[0,\bar \rho)$, $\bar \rho \le 1$. Specifically, for each $v>v_{max}$ the $n$-th system is stable (for all sufficiently large $n$) and the SSAI holds, with the limiting distribution concentrated on the ML-MFP $x^{*,v}$ (which is the unique 
left-regulated proper ML-FP for this speed), and load 
$\rho(v) =x^{*,v}_0 \uparrow \bar \rho$ as $v \downarrow v_{\max}$.
(Note that if, in addition, Condition~\ref{cond-unique} holds for all $v$, then, by Theorem~\ref{thm-free-wave-to-regulated-unique}(v), 
$\bar \rho =1$, and therefore SSAI-FRL holds.)
\end{thm}

{\em Proof of Theorem~\ref{thm-free-wave-to-ssai}.} 
By Theorem~\ref{thm-free-wave-to-regulated-unique}, there is the family of unique, left-regulated ML-FP (=ML-MFP) $x^{*,v}$,
with load $\rho(v)=x^{*,v}_0$ continuously strictly increasing from $0$ to $\bar \rho$ as $v$ decreases from $+\infty$ to $v_{max}$;
and Condition~\ref{cond-unique} holds for any $v>v_{max}$. Denote by $v(\rho)$ the inverse of $\rho(v)$.

Pick any value of load $\rho < \bar \rho$. Note that $v(\rho) > v_{max}$.
We can and do choose a sequence of systems, in $n$, 
such that they are stable, 
 $v^n \to v$, 
and load $\rho^n = \E x^n_0(\infty) \to \rho$. Condition $v\ge v(\rho)$ (and then $v > v_{max}$) must necessarily hold --
otherwise we would have that $x^{*,v}$ strictly dominates $x^{*,v(\rho)}$, and then
$\lim \E x^n_0(\infty) \ge \rho(v) > \rho$. 
Therefore, Condition~\ref{cond-unique} holds for this $v$. Then, by Proposition~\ref{prop-ssai-key},
SSAI holds, and $v$ must be equal to $v(\rho)$. 
$\Box$

\begin{cor}
\label{thm-free-wave-to-ssai-cor1}
Suppose Assumptions~\ref{cond-second-moment}, \ref{cond-ntriv-add-ii} and \ref{cond-ntriv-add-i} hold. 

(i) If $v_{min} = v_{max}=v_*$, then
$\bar \rho$ in Theorem~\ref{thm-free-wave-to-ssai} is the supremum of the loads $x^{*,v}_0$ of all 
proper left-regulated ML-MFP $x^{*,v}$ in the frame $[0,\infty)$, under different speeds $v$.

(ii) If Assumption~\ref{assump-101} holds, then conclusions of (i) hold. (This result was obtained as theorem 10 in \cite{St2021-coc}.)
\end{cor}

{\em Proof.} 
(i) If $v_{min} = v_{max}=v_*$, then a proper left-regulated ML-FP does not exist for $v<v_*$. 
And we know that a unique proper left-regulated ML-FP ( = ML-MFP) does exist for $v > v_{max}=v_*$.

(ii) If Assumption~\ref{assump-101} holds, then $v_{min} = v_{max}=v_*$.
$\Box$

\begin{cor}
\label{thm-free-wave-to-ssai-cor2}
Suppose Assumptions~\ref{cond-second-moment}, \ref{cond-ntriv-add-ii} and \ref{cond-ntriv-add-i} hold. If, in addition,  D-monotonicity holds (and then $v_{min} = v_{max}=v_*$) we have: $\bar \rho$ in Theorem~\ref{thm-free-wave-to-ssai} is equal to $1$, so that SSAI-FRL holds. 
(This result also follows from theorems 4 and 6 in \cite{St2021-coc}.)
Furthermore, as $v \downarrow v_*$,  the (appropriately shifted version of) ML-MFP $x^{*,v}$ converges to the unique free ML-FP for speed $v_*$.
\end{cor}

{\em Proof.} Under D-monotonicity, Condition~\ref{cond-unique} holds for all $v$. 
 $\Box$
 
\begin{thm}
\label{thm-ssai-right}
Suppose Assumptions~\ref{cond-second-moment}, \ref{cond-ntriv-add-ii} and \ref{cond-ntriv-add-i} hold.
If free ML-FP speed range is $[v_{min},v_{max}]$, then SSAI  
holds for the right-regulated systems with frame $(-\infty,0]$ for speeds $v< v_{min}$.
Specifically, for each $v< v_{min}$ the $n$-th system is stable (for all sufficiently large $n$) and the SSAI holds, with the limiting distribution concentrated on the ML-MFP $x^{*,v}$ (which is the unique 
right-regulated proper ML-FP for this speed).
Furthermore, as $v \uparrow v_{min}$,  (a shifted version of) $x^{*,v}$ converges to the unique proper free ML-FP for speed $v_{min}$.
\end{thm}

We will not give a full proof of Theorem~\ref{thm-ssai-right} in this paper, because it is analogous to the proof of Theorem~\ref{thm-free-wave-to-ssai} 
-- we only give a sketch in Section~\ref{sec-right-regulated-ssai-proof}. But we want to emphasize that Theorem~\ref{thm-ssai-right} for the right-regulated system gives SSAI property which is, in a sense, stronger than that the SSAI property in Theorem~\ref{thm-free-wave-to-ssai} for the left-regulated system.
Specifically, {\em without any extra conditions},
Theorem~\ref{thm-ssai-right} shows that, as $v \uparrow v_{min}$, the ML-MFP $x^{*,v}$ converges to the proper free ML-FP for speed $v_{min}$.
In contrast, Theorem~\ref{thm-free-wave-to-ssai} for the left-regulated system shows that, 
as $v \downarrow v_{max}$, the ML-MFP $x^{*,v}$ converges to the proper free ML-FP for speed $v_{max}$, {\em only when additional
Condition~\ref{cond-unique} holds for all $v$} 
(for which D-monotonicity is sufficient.)

\begin{thm}
\label{thm-vn-limit}
Consider a free system. Suppose Assumptions~\ref{cond-second-moment}, \ref{cond-ntriv-add-ii} and \ref{cond-ntriv-add-i} hold.

(i) Suppose the centered free process $\mathring x^n(\cdot)$ is positive recurrent, and therefore 
the average velocity $v_n$ of unanchored particle system is well-defined -- see \eqn{eq-vn-def-wp1}-\eqn{eq-vn-def-stat}.
(Positive recurrence holds, for example, if, in addition, Assumption~\ref{assump-small} holds -- see Theorem~\ref{th-closeness-free}.)
Then, $\liminf_n v_n \ge v_{min}$ and $\limsup_n v_n \le v_{max}$, where $[v_{min},v_{max}]$ is the 
free proper ML-FP speed range.

(ii) More generally, denote\\
$v_{min,n} = \sup\{v ~|~ \mbox{$\forall \nu \in (0,1), ~\liminf_{t\to\infty} \Delta^n(t;\nu)/t\ge v$, w.p.1} \}$, \\
$v_{max,n} = \inf\{v ~|~ \mbox{$\forall \nu \in (0,1), ~\limsup_{t\to\infty} \Delta^n(t;\nu)/t\le v$, w.p.1} \}$. \\
Then,
$\liminf_n v_{min,n} \ge v_{min}$  
and $\limsup_n v_{max,n} \le v_{max}$.
\end{thm}

{\em Proof of Theorem~\ref{thm-vn-limit}.}  
(i) We first  prove that $\liminf_n v_n \ge v_{min}$. By Theorem~\ref{thm-free-wave-spectrum}(iii), 
 for any speed $v < v_{min}$ there exists a right-regulated
ML-FP $\tilde x$. Let us choose $v$ slightly smaller than $v_{min}$. If we simultaneously scale down $v$ and the arrival rates $\{\sigma_j\}$ by a factor $\gamma$ slightly smaller than $1$, we obtain speed $\tilde v = \gamma v$, which is still only slightly smaller than $v_{min}$, and the arrival rates
$\{\tilde \sigma_j\} =  \gamma \{\sigma_j\}$, for which $\tilde x$ is still an ML-FP. Now, just like in the proof of 
Theorem~\ref{thm-free-wave-speed-as-lower-bound},  let us fix a small $\delta>0$
and consider unanchored ML $x^\ell$, starting from state $x^\ell(0)=(1-\delta) x^\emptyset$. 
The velocity of $x^\ell$, restricted to $\nu$-quantiles with $\nu \in (\delta,1)$, is equal to the velocity of
$(0,\delta)$-reduced system; the arrival rates in the $(0,\delta)$-reduced system are only slightly different from the original
rates $\{\sigma_j\}$, and therefore dominate $\{\tilde \sigma_j\}$, if we choose $\delta$ small enough. 
We see that the velocity of $x^\ell$, restricted to $\nu$-quantiles with $\nu \in (\delta,1)$, is at least $\tilde v$.
And recall that $(\tilde v, \delta)$ can be arbitrarily close to $(v_{min},0)$.

Consider now the $n$-th pre-limit unanchored process $x^n(\cdot)$ (with $n$ particles) in the stationary regime,
in the sense that the centered process $\mathring x^n(\cdot)$ is in the stationary regime.
WLOG, let the $\delta$-quantile of the process (random) initial 
 state $x^n(0)$ be located exactly at $0$.
Consider $\nu$-quantile of $x^n(0)$, with $\nu \in (\delta,1)$.
 The $\nu$-quantile's expected distance from the $\delta$-quantile is uniformly bounded in $n$. (Here we use the fact
that, by Theorem~\ref{th-closeness-free}, $\E\Phi_1(\mathring x^n(\infty))$ is finite, uniformly bounded.) If we consider a different initial state,
with $\delta$-fraction of particles moved back to $-\infty$ and the remaining $(1-\delta)$-fraction of particles moved back to $\delta$-quantile,
 then, by monotonicity,
we obtain a lower bounding process which is dominated by $x^n(t)$ at all times $t$. But, as $n\to\infty$, the ML for this lower bounding process is exactly
$x^\ell(\cdot)$, whose velocity (restricted to $\nu$-quantiles with $\nu \in (\delta,1)$)
 is at least $\tilde v$.
Using this, \eqn{eq-vn-def-stat} and \eqn{eq-vn-def-wp1}, we can write
$$
\liminf_n v_n = \liminf_n \E [\zeta^n(t;\nu) - \zeta^n(0;\nu)]/t \ge \tilde v + o(1),
$$ 
where $\zeta^n(t;\nu)$ is the location of $\nu$-quantile at time $t$, and $o(1)$ is a function of $t$, vanishing as $t\to\infty$.
(The $o(1)$ term accounts for the quantile $\zeta^n(0;\nu)$ being ``ahead'' of the $\nu$-quantile of the lower bounding process, which 
is located at $\zeta^n(0;\delta) = 0$. And, as already noted above, $\E [\zeta^n(0;\nu) - \zeta^n(0;\delta)]$ is uniformly bounded in $n$.)
Therefore, $\liminf_n v_n \ge \tilde v$. And since $\tilde v$ can be arbitrarily close to $v_{min}$, we conclude that $\liminf_n v_n \ge v_{min}$.

The proof of $\limsup_n v_n \le v_{max}$ is by contradiction. Suppose there is a subsequence of $n$, along which $\lim v_n=v' > v_{max}$.
Consider the corresponding subsequence of left-regulated anchored systems with speed $v \in (v_{max},v')$. This system must be unstable for large $n$,
because $v_n > v$, and average velocity of particles in the left-regulated system dominates that in the free system.
On the other hand, for large $n$ the left-regulated system has to be stable by Theorem~\ref{thm-free-wave-to-ssai} -- a contradiction.
 
(ii) When $v_n$ is not necessarily well-defined, we
can introduce upper 
and lower bounding systems, satisfying Assumption~\ref{assump-small}, and with $v_{min}$ and $v_{max}$ being, respectively, slightly smaller and larger than in original system. 
We can then use the special result in (i) for both, and then use Theorem~\ref{thm-free-wave-vmin-vmax-monotone}. The details are as follows.

The construction of an upper bounding system is straightforward. 
For example, we add a new class, with very small arrival rate, and i.i.d. exponential component sizes.

To construct a lower bounding system, consider a class, satisfying Assumption~\ref{cond-ntriv-add-ii}. Then, WLOG, there is a class $j$ such that
there are $m_j$ components, $0 \le m_j < k_j$, which are exactly $0$, and the remaining $d_j - m_j$ components have random sizes which are
between $\epsilon>0$ and $\epsilon_2>\epsilon$. We define a lower bounding system by replacing class $j$ by two classes as follows: the first one occurs with high probability $1-\delta_1 < 1$ and has the same component size distribution $F_j$ as  class $j$; the second one 
occurs with small probability $\delta_1>0$, and the component size distribution is such that we, first, take a realization of sizes according 
to $F_j$ and then replace each non-zero component size realization $s$ by realization of the uniform distribution 
in $[0,s]$. The second of the new classes (of the lower bounding system) satisfies conditions of Lemma~\ref{lem-small-sufficient}. Therefore, the lower bounding system satisfies Assumption~\ref{assump-small}. 
Now, as $\delta_1 \downarrow 0$ along some sequence, conditions of Theorem~\ref{thm-free-wave-vmin-vmax-monotone} are satisfied,
and so we can choose a lower bounding system with $v_{min}$ being arbitrarily close to that of the original system. 
$\Box$

\begin{thm}
\label{thm-vn-iid-ihr}
Consider a free system. Suppose for each job class $j$ the distribution $F_j$ is such that the components are i.i.d. with IHR
or, more generally, $F_j$ is a mixture of such distributions. Suppose that $k_j < d_j$ for at least one class $j$, and 
Assumption~\ref{cond-second-moment} holds.
Then we have:

(i) $\lim_n v_{min,n} = \lim_n v_{max,n} = v_*$, where $v_*=v_{min}=v_{max}$ is the speed of the unique free ML-FP. 

(ii) Suppose, in addition, that $v_n$ is well-defined, i.e., the free process $\mathring x^n(\cdot)$ is positive recurrent. 
(Positive recurrence holds, for example, if, in addition, Assumption~\ref{assump-small} holds -- see Theorem~\ref{th-closeness-free}.) 
Then, 
$\lim_n v_n  = v_*$.
\end{thm}

Note that the existence of class $j$ with $k_j < d_j$, and the (mixture of) i.i.d. IHR assumption for 
this class, imply Assumption~\ref{cond-ntriv-add-ii}. So, conditions of Theorem~\ref{thm-vn-iid-ihr}
imply that both Assumptions~\ref{cond-second-moment} and \ref{cond-ntriv-add-ii} hold.
We also remark that each distribution $F_j$ satisfies the i.i.d. IHR assumption, or is a mixture of a {\em finite number} of such distributions,
then Assumption~\ref{cond-second-moment} holds automatically. However, if at least one distribution $F_j$ is a mixture of 
more than a finite number of i.i.d. IHR distributions, Assumption~\ref{cond-second-moment} in Theorem~\ref{thm-vn-iid-ihr}
needs to be made explicitly.

{\em Proof of Theorem~\ref{thm-vn-iid-ihr}.} As noted earlier, the theorem conditions imply that both
Assumptions~\ref{cond-second-moment} and \ref{cond-ntriv-add-ii} hold. Then  ML-FP uniqueness 
holds by Theorem~\ref{thm-free-wave-unique-cond}(ii).

(i) We would like to apply Theorem~\ref{thm-vn-limit}(ii). However, we cannot do it directly, because Assumption~\ref{cond-ntriv-add-i}
does not necessarily hold. To get around this difficulty, 
we can construct upper and lower bounding systems, satisfying in addition  Assumption~\ref{cond-ntriv-add-i}, with $v_{max}$ and $v_{min}$
being respectively slightly larger and slightly smaller than $v_*$. (This is done the same way as in the proof of Theorem~\ref{thm-vn-limit}(ii). 
But here the purpose of the construction is not to make the process positive recurrent by making it satisfy Assumption~\ref{assump-small} --
it is instead to make the process satisfy Assumption~\ref{cond-ntriv-add-i}. Same construction does that.)
Then we can apply Theorem~\ref{thm-vn-limit}(ii) to those lower and upper bounding systems, which then leads to
$\lim_n v_{min,n} = \lim_n v_{max,n} = v_*$.

(ii) If $v_n$ is well-defined, then (i) implies $\lim_n v_n  = v_*$. $\Box$

\section{Proofs}
\label{sec-proofs-new}

\subsection{Proof of Proposition~\ref{prop-ml-result}}
\label{sec-ml-proof}

The proof of the convergence to a unique ML repeats the proof of theorem 7 in \cite{St2021-coc},
except here we consider a more general system, which allows regulation on either side, or on both sides, or no regulation at all.
And also Proposition~\ref{prop-ml-result} is formulated in terms of unanchored system, instead of corresponding anchored system.
The fact that this ML must be an MM follows from the interpretation of ML given immediately after the proof of theorem 7 in \cite{St2021-coc}.
For the uniqueness of MM, consider any two MM $x(\cdot)$ and $y(\cdot)$ with the same initial state $x(0)=y(0)$, and 
note that
from \eqn{eq-h-Lipschitz} we easily obtain 
\beql{eq-deriv-universal-2}
 \frac{d}{d t} \left\| x(t) - y(t) \right\| \le \lambda C \| x(t) - y(t) \|,
\eeql
as long as $\| x(t) - x(t) \| >0$. But then, by Gronwall inequality, $\| x(t) - y(t)\| \equiv 0$.
This proves the uniqueness.
$\Box$

\subsection{Proof of Lemma~\ref{lem-hor-dist-no-local-max-new}}
\label{sec-monotone-diff}

 The proof is by contradiction. If lemma does not hold, then the following holds for some $\zeta$ satisfying $\zeta \le \chi$ and
 $\zeta \in (b,1]$, 
 and some $\beta \in (b,\zeta)$:
 \beql{eq-local-max1-new}
 \psi_\zeta = \max_{u\in [\zeta,1]} \psi_u
 \eeql
 and 
 \beql{eq-local-max2-new}
 \psi_u < \psi_\zeta, ~~\forall u\in (\beta,\zeta).
 \eeql
 For the rest of the proof let us shift $x$, along with its frame, to the right by $\psi_\zeta$. (This means the shift is actually to the left if $\psi_\zeta<0$.)
 With this ``new'' $x$:
 $\psi_\zeta = 0$, $\psi_u \le 0$ for $u \ge \zeta$, and $\psi_u < 0$ for $u \in (\beta,\zeta)$; consequently, denoting
 $\hat w = \hat \phi_\zeta$ and $\Delta(y)=x_{\hat w + y} - \hat x_{\hat w + y}$, we have
 that $\Delta(y) \ge 0$ for $y \le 0$, $\Delta(0)=0$, and $\Delta(y) > 0$ in some small interval $y \in (0, \gamma]$.
 For $y \in [0, \gamma]$, denote
 $$
 D(y) = \max_{0\le \xi \le y} \Delta(\xi),
 $$
 which is a non-decreasing Lipschitz function with $D(0)$ and $D(y)>0$ for $y>0$. Since $\Delta(y)$ is also Lipschitz,
 we have the following properties: the derivative $D'(y)$ exists a.e., and a.e. $D'(y)>0$ implies that $D(y)=\Delta(y)$ and $D'(y)=\Delta'(y)$. 
 Consider a point $y>0$ such that $D(y)=\Delta(y)$ and $D'(y)=\Delta'(y)>0$. 
 Using the fact that $v \le \hat v$ and 
 $\{\sigma_j\} \ge \{\hat \sigma_j\}$, 
 and the properties of function $h(x_{(-\infty,w]})$, in particular \eqn{eq-h-Lipschitz}
 and monotonicity,
 we obtain that 
 \beql{eq-delta-ineq}
 \Delta'(y) \le \frac{\hat \lambda \bar d^2}{\hat v} \Delta(y), ~~\hat \lambda = \sum_j \hat \sigma_j.
 \eeql
 (Here we compare $\hat x_w$ to the function $\tilde x_w$, which coincides with $\hat x_w$ for $w \le \hat w$,
 coincides with $x_w$ for $\hat w < w < \hat w + y$, and with $\tilde x_{\hat w +y} = \hat x_{\hat w +y}$,
 so that it possibly makes a jump down at $\hat w +y$, with the jump size being $\Delta(y)$.
 This jump down is what accounts for the RHS of \eqn{eq-delta-ineq}, by the same argument as in lemmas 12-13 in \cite{St2021-coc}.)
 From \eqn{eq-delta-ineq} we see that a.e. 
 $$
 D'(y) \le \frac{\hat \lambda \bar d^2}{\hat v} D(y).
 $$
 But then, by Gronwall inequality, $D(y)=0$. The contradiction completes the proof. $\Box$

\subsection{Proof of Theorem~\ref{thm-bounded-first-moment}}
\label{sec-closeness-proof}

If $x^n(t)$ is the state at time $t$, it can be equivalently described as
$p^n(t)=(w_1(t), w_2(t), \ldots, w_n(t); z(t))$, where $z(t) = - \bar x^n(t) \le 0$ is the location of the left boundary w.r.t. to the mean $\bar x^n(t)$,
and 
$w_1, w_2, \ldots, w_n \ge z$ are the locations of the $n$ particles (not necessarily ordered in any way)
w.r.t. the mean $\bar x^n(t)$. 
We will slightly abuse notation by defining $\Phi_\ell(p^n(t))$ as $\Phi_\ell(\mathring x^n(t))$ for the corresponding
$\mathring x^n(t)$, that is
$$
\Phi_\ell(p^n(t)) = \frac{1}{n} \sum_{i=1}^n |w_i(t)|^\ell.
$$

Since a stationary distribution is invariant w.r.t. simultaneous scaling of speed $v$ and $\lambda^n$,
let us assume that $v=1$ WLOG.
Then the evolution of $p^n(t)$ is as follows. Note that here the particle locations $w_i(t)$ may be negative as well as non-negative.
Point $z(t)$ serves as a regulation boundary, which evolves in time; $z(t)\le 0$ at all times.
Denote by $I_z(t)$ the (possibly empty) subset of particles, 
which are located exactly at the regulation boundary  $z(t)$ at time $t$.
Between the times of the particle jumps (job arrivals), 
the  boundary $z(t)$ moves right at the constant non-negative speed $(n-|I_z(t)|)/n$; 
the particles in $I_z(t)$ (those located exactly at the boundary $z(t)$) 
stay at  the boundary and therefore move with it to the right at constant non-negative speed $(n-|I_z(t)|)/n$; 
each particle that is not at $z(t)$ (i.e. not within subset $I_z(t)$) moves to the left at the constant non-positive speed $-|I_z(t)|/n$;
when a particle ``hits'' boundary $z(t)$, it joins the set $I_z(t)$. (The average of the particle locations stays at $0$, as it should by
the $p(t)$ definition.)
It is easy to see that, for any $\ell\ge 1$, $\Phi_\ell(p^n(t))$ is absolutely continuous with the derivative (existing almost everywhere) 
$(d/dt)\Phi_\ell(p^n(t)) \le 0$, and moreover $(d/dt)\Phi_\ell(p(t)) < 0$ as long as  $0< |I_z(t)| < n$. At a time when one or more particles jump (upon a job arrival), the following occurs. Let $\kappa_i(t)$ be the amount of new workload added to server $i$ by the
job arriving at time $t$ (which may be non-zero only for the servers within the selection set of the job). Then, particle $i$ jump size at $t$ 
is $\zeta_i=\kappa_i(t) - \sum_s \kappa_s(t)/n$ (which may be positive or negative); the point $z(t)$ jumps (left) by 
$ - \sum_s \kappa_s(t)/n$.

Consider the following function $G^{(n)}(\mathring x^n), ~\mathring x^n \in \mathring\cx^{(n)},$ 
$$
G^{(n)} (\mathring x^n) = G^{(n)} (p^n)
= 2 \lambda^n n \frac{1}{n} \sum_i w_i \E  \zeta_i^{(\mathring x^n)},
$$
where $\zeta_i^{(\mathring x^n)}$ is the random jump size (which can have any sign) of particle $i$ upon a job arrival when the state is 
$\mathring x^n$. (The sizes $\zeta_i^{(\mathring x^n)}$ are dependent across $i$, of course.) 
As we will see, function $G^{(n)} (p^n)$ can be thought of as the ``first-order approximation of the generator of process $p^n(\cdot)$, applied to function $\Phi_2(p^n)$;'' but we do not even claim that $\Phi_2(p^n)$ is within the generator domain.

For future reference, note that, for each $n$, $G^{(n)} (\mathring x^n)$ is continuous in $\mathring x^n \in \mathring\cx^{(n)}$.
Also note that each 
$\E  \zeta_i^{(\mathring x^n)}$ is the quantity of the order $O(1/n)$, which motivates the definition
\beql{eq-bar-zeta-def}
\bar \zeta_i^{(\mathring x^n)} \doteq n \E \zeta_i^{(\mathring x^n)}.
\eeql
We can then write
$$
G^{(n)} (\mathring x^n) = G^{(n)} (p)
= 2 \lambda^n \frac{1}{n} \sum_i w_i \bar \zeta_i^{(\mathring x^n)}.
$$
We define the function $\bar\zeta^{(\mathring x^n)}(w), ~w\in \R,$ as follows: 
$\bar \zeta^{(\mathring x^n)}(w) = \bar \zeta^{(\mathring x^n)}_i$, where $i$ is the particle whose location $w_i$ is the closest to $w$ on the left; we also adopt a conventions that, if $w_i$ is the location 
of the left-most particle, then $\bar\zeta^{(\mathring x^n)}(w) = \bar \zeta_i^{(\mathring x^n)}$ for all $w< w_i$.
Clearly, function $\bar\zeta^{(\mathring x^n)}(w)$ is a piece-wise constant non-increasing function, and we can write
\beql{eq-gen1}
G^{(n)} (\mathring x^n) = G^{(n)} (p^n)
= 2 \lambda^n \int_{-\infty}^{\infty} \bar\zeta^{(\mathring x^n)}(w) w d[-\mathring x^n_w] \le 0,
\eeql
where the inequality in \eqn{eq-gen1} follows because
$$
\int_{-\infty}^{\infty} \bar\zeta^{(\mathring x^n)}(w) w d[-\mathring x^n_w] 
= \int_{-\infty}^{\infty} [\bar\zeta^{(\mathring x^n)}(w) - \bar\zeta^{(\mathring x^n)}(0)] w d[-\mathring x^n_w] 
+ \bar\zeta^{(\mathring x^n)}(0) \int_{-\infty}^{\infty}  w d[-\mathring x^n_w],
$$
$\bar\zeta^{(\mathring x^n)}(w) - \bar\zeta^{(\mathring x^n)}(0)$ is non-increasing with value $0$ at $0$,
and $\mathring x^n \in \mathring{\cx}_1$
(and then $\int_{-\infty}^{\infty} w d[-\mathring x^n_w] = 0$).

Next, we claim the following property: there exists a sufficiently large $C>0$ and some $\epsilon > 0$, such that,
uniformly in all sufficiently large $n$ and all $\mathring x^n \in \mathring\cx^{(n)}$ with $\Phi_1(\mathring x^n) \ge C$,
\beql{eq-drift}
G^{(n)} (\mathring x^n) = G^{(n)} (p^n) \le -\epsilon \Phi_1(\mathring x^n).
\eeql
The proof of \eqn{eq-drift} is given in Section~\ref{sec-eq-drift}.

From \eqn{eq-drift} and \eqn{eq-gen1} we obtain that, uniformly in all sufficiently large $n$,
\beql{eq-drift2}
G^{(n)}(\mathring x^n) = G^{(n)}(p^n) \le -\epsilon \Phi_1(\mathring x^n) +\epsilon C.
\eeql

Denote by $\Phi_2^{(C_2)}(p^n) = \Phi_2(p^n) \wedge C_2$ the function $\Phi_2$ truncated at level $C_2$.
This is a continuous bounded function. However, we cannot claim that $\Phi_2^{(C_2)}(\cdot)$ is within 
the domain 
of the generator of process $p^n(\cdot)$, because state $p^n(t)$ may change continuously between
job arrivals (when/if one or more particles are located at the regulation boundary), and 
$\Phi_2^{(C_2)}(\cdot)$ is not differentiable. To deal with this situation, we will work with the ``generator upper bound''
$A^{(n)}_{C_2}$, given in the following

\begin{lem}
\label{lem-generator-bound}
Let $C_2 >0$ and $n$ be fixed.
Let $\Delta^{(n)}_{C_2}(p^n)$ denote the expected change of $\Phi_2^{(C_2)}(p^n)$ when a job arrival occurs at state $p^n$,
$A^{(n)}\Phi_2^{(C_2)}(p^n) \doteq \lambda^n n \Delta^{(n)}_{C_2}(p^n)$,
$$
    A^{(n)}_{C_2}\Phi_2^{(C_2)}(p^n) \doteq
    \begin{cases}
       A^{(n)}\Phi_2^{(C_2)}(p^n), & \text{if } \Phi_2 (p^n) \le C_2 \\
      0, & \text{otherwise}
    \end{cases}
$$
    and
$$
B^t \Phi_2^{(C_2)}(p^n) \doteq \frac{P^t - I}{t} \Phi_2^{(C_2)}(p^n).
$$
Then
$$
\limsup_{t \downarrow 0} \sup_{p^n} [B^t \Phi_2^{(C_2)}(p^n) - A^{(n)}_{C_2}\Phi_2^{(C_2)} (p^n)] \le 0.
$$
\end{lem} 
The proof of Lemma~\ref{lem-generator-bound} is in Section~\ref{sec-generator-bound2}.

Next, we claim the following fact:
there exists $C_1>0$ such that for any fixed $C_2>0$, 
uniformly in all large $n$ and $p^n$ such that $\Phi_2(p^n) \le C_2$, we have
\beql{eq-generator-bound}
A^{(n)}_{(C_2)} \Phi_2^{(C_2)}(p^n) \le G^{(n)}(p^n) + C_1,
\eeql
and then
\beql{eq-generator-bound5}
A^{(n)}_{(C_2)} \Phi_2^{(C_2)}(p^n) \le  -\epsilon \Phi_1(p^n) + C_3,
\eeql
with $C_3 = \epsilon C + C_1$.
The proof of \eqn{eq-generator-bound} is given in Section~\ref{sec-generator-bound}.

Bound \eqn{eq-generator-bound5} in turn implies that
 for any fixed $C_2>0$, 
\beql{eq-generator-bound2}
A^{(n)}_{(C_2)} \Phi_2^{(C_2)}(p^n)  \le [-\epsilon \Phi_1(p^n) + C_3] \bI \{\Phi_2(p^n) \le C_2\}.
\eeql

Recalling that $p^n(\infty)$ is the random value of $p^n(t)$ in the stationary regime, 
we have for all large $n$, and using Lemma~\ref{lem-generator-bound}, we can write:
$$
0 = \lim_{t\downarrow 0} \E B^t \Phi_2^{(C_2)}(p^n(\infty)) \le 
\E A^{(n)}_{(C_2)}  \Phi_2^{(C_2)}(p^n(\infty)) \le 
\E \left[ (-\epsilon \Phi_1(\mathring x^n(\infty)) + C_3) \bI \{\Phi_2(\mathring x^n(\infty)) \le C_2\} \right]
$$
$$
\le -\epsilon \E [\Phi_1(\mathring x^n(\infty))  \bI \{\Phi_2(\mathring x^n(\infty)) \le C_2\}] + C_3,
$$
and then
$$
\E [\Phi_1(\mathring x^n(\infty))  \bI \{\Phi_2(\mathring x^n(\infty))  \le C_2\}] \le C_3/\epsilon.
$$
Letting $C_2 \uparrow \infty$, we finally obtain that
$$
\E \Phi_1(\mathring x^n(\infty)) \le C_3/\epsilon
$$
for all sufficiently large $n$, and then 
$$
\E \Phi_1(\mathring x^n(\infty)) \le \bar C
$$
holds for all $n$ for some large $\bar C>0$.
$\Box$

\subsubsection{Proof of \eqn{eq-drift}.}
\label{sec-eq-drift}

The definition of $\bar \zeta_i=\bar \zeta_i^{(\mathring x^n)}$ in \eqn{eq-bar-zeta-def} can be interpreted as follows:  
$\bar \zeta_i$ is the expected actual jump size 
$\E [\kappa_i | i\in S]$ of particle $i$, {\em conditioned on it being selected by a job}, multiplied by constant $[n\P\{i\in S\}]$
(which is close to the constant $\sum_j \pi_j d_j$ for all large $n$), and then centered by the total expected particles' displacement
 $\E \sum_s \kappa_s$ due to a job arrival. 

The proof is by contradiction. Suppose property \eqn{eq-drift} does not hold. Then, we can and do choose a subsequence of $n\to\infty$, and corresponding $\mathring x^n$, so that 
along this subsequence $c_n =\Phi_1(\mathring x^n) \uparrow \infty$ and 
\beql{eq-contr}
G^{(n)} (\mathring x^n) / c_n \to 0.
\eeql
Using the continuity of $G^{(n)} (\mathring x^n)$ in $\mathring x^n \in \mathring\cx^{(n)}$, we also can and
 do choose this subsequence so that all particle locations $w_i$ are distinct (i.e., $\mathring x^n$ has exactly $n$ jumps,
at points $w_i$, each jump is by $1/n$). 
For a fixed $\nu \in (0,1/4)$ and each $n$, denote 
$$
a^+_{\nu,n}= a^+_{\nu,n}(\mathring x^n)= \max \left\{w^+ \ge 0 ~|~ \int_{w^+}^\infty w d[-\mathring x^n_{w}]  \ge \nu c_n \right\},
$$
$$
a^-_{\nu,n}= a^-_{\nu,n}(\mathring x^n)= \min \left\{w^- \le 0 ~|~ \int_{-\infty}^{w^-} |w| d[-\mathring x^n_{w}]  \ge \nu c_n \right\}.
$$
It is easy to see that for a sufficiently small fixed $\delta>0$, uniformly in $n$ and $\nu \in (0,1/4)$,
$$
a^+_{\nu,n} \ge \delta c_n, ~~ a^-_{\nu,n} \le -\delta c_n.
$$
(Otherwise, $\int_{-\infty}^{0} |w| d[-\mathring x^n_{w}] = \int_0^{\infty} w d[-\mathring x^n_{w}] = c_n/2$ could not hold.)
Denote 
$$
q^+ = \mathring x^n_{a^+_{\nu,n}}, ~~ q^- = 1-\mathring x^n_{a^-_{\nu,n}}.
$$
It is easy to see that for any $\delta_1>0$ we can fix a sufficiently small $\nu>0$ such that, uniformly in all large $n$, 
$$
q^+ \le \delta_1, ~~ q^- \le \delta_1.
$$
(If not, since both $a^+_{\nu,n}/c_n$ and $a^-_{\nu,n}/c_n$ are bounded away from $0$,
$\nu$ could not be small.) Let us fix any $\delta_1 \in (0,1/2)$, and then choose a corresponding $\nu$.

Observe that the following holds for each large $n$. The distance between the particles located
at $a^-_{\nu,n}$ and $a^+_{\nu,n}$, let us label them $i^-$ and $i^+$, 
 is at least $2\delta c_n$, and there is at least $(1-2\delta_1)n/2$ (i.e. a non-zero fraction) of particles, located strictly between particles
 $i^-$ and $i^+$. 
 
 Using coupling and Assumptions~\ref{cond-ntriv-add-ii} and \ref{cond-ntriv-add-i} (this is the only place in the proof where these assumptions is used),
 we obtain that there exists $\delta_2>0$ such that
 \beql{eq-drift-diff}
 \bar\zeta^{(\mathring x^n)}_{i^-} - \bar\zeta^{(\mathring x^n)}_{i^+} \ge \delta_2.
 \eeql
 Indeed, when a job arrival selects both $i^-$ and $i^+$, then, conditioned on this event, 
 by monotonicity, $i^-$ receives at least as much average new workload as $i^+$.
 Also, clearly, the probability that both $i^-$ and $i^+$ are selected by an arriving job, under the condition that $i^-$
 (respectively, $i^+$) is selected vanishes as $n\to\infty$.
 Thus, to show \eqn{eq-drift-diff}, is suffices to consider two probability spaces, $E^-$ and $E^+$, 
 obtained by conditioning on the two non-intersecting events, when, respectively, $i^-$ (but not $i^+$) and $i^-$ (but not $i^+$) is selected
 by a job.
We can couple the constructions of the spaces $E^-$ and $E^+$, so that, w.p.1, 
$i^-$ (on $E^-$) and $i^+$ (on $E^+$) receive equal component sizes, while 
all other selected servers and their component sizes are equal.
 At least one of the following two conditions holds for $\mathring x^n$: (a) at least $[(1-2\delta_1)n/2]/2$
 particles are located in $[0,a^+_{\nu,n})$ or (b) at least $[(1-2\delta_1)n/2]/2$
 particles are located in $(a^-_{\nu,n},0]$. In the case (a), 
 consider an arrival of a job of class $j$, satisfying Assumption~\ref{cond-ntriv-add-i}.
 Then, with probability bounded away from $0$ 
 all other selected particles are in $[0,a^+_{\nu,n})$; and then, Assumption~\ref{cond-ntriv-add-i} easily implies \eqn{eq-drift-diff},
 because, with bounded away from $0$ probability, $i^-$ will receive new workload equal to its component size, while $i^+$ will receive 
 a strictly smaller (by a bounded away from $0$ quantity) new workload. 
 In the case (b), consider an arrival of a job of class $j$, satisfying Assumption~\ref{cond-ntriv-add-ii}.
 With probability bounded away from $0$ 
 all other selected particles are in $(a^-_{\nu,n},0]$; and then, Assumption~\ref{cond-ntriv-add-ii} easily implies \eqn{eq-drift-diff},
 because, with bounded away from $0$ probability, $i^-$ will make some positive (bounded away from $0$) jump forward,
 while $i^+$ will not jump. This completes the proof of \eqn{eq-drift-diff}.
 
 From \eqn{eq-drift-diff} we see that 
 $$
 \mbox{either}~ \bar\zeta^{(\mathring x^n)}(a^-_{\nu,n}) - \bar\zeta^{(\mathring x^n)}(0) \ge \delta_2/2~~\mbox{or}~~
 \bar\zeta^{(\mathring x^n)}(a^+_{\nu,n}) - \bar\zeta^{(\mathring x^n)}(0) \le -\delta_2/2.
 $$
 In either case, $G^{(n)} (\mathring x^n) \ge (\delta_2/2) \nu c_n$, which contradicts \eqn{eq-contr}.
$\Box$

\subsubsection{Proof of Lemma~\ref{lem-generator-bound}.}
\label{sec-generator-bound2}

Function $\Delta^{(n)}_{C_2}(p^n)$, and then 
$A^{(n)} \Phi_2^{(C_2)}(p^n)$,  is continuous by Lemma~\ref{lem-monotone-cont-basic},
and $A^{(n)} \Phi_2^{(C_2)}(p^n) \le 0$ when $\Phi_2 (p^n) > C_2$. 
Consider the process evolution over a small interval $[0,t/n]$, with $p^n(0)=p^n$ such that $\Phi_2(p^n) \le C_2$.
Note that, if we are interested in a limit of $B^{t/n} \Phi_2^{(C_2)}(p^n)$ as $t\downarrow 0$,
we can ``ignore'' the event when more than one job arrival occurs -- its contribution into the expected change of  $\Phi_2^{(C_2)}$
is $o(t)$. (Because jobs arrive as Poisson process of rate $\lambda^n n$, and $\Phi_2^{(C_2)}$ is bounded.) 
Note that as $p^n$ changes continuously due to, maybe, some particles being on the boundary, the value of
$\Phi_2(p^n)$ cannot increase. (Movement of the particles can be decomposed into pairs of particles -- on and off the boundary --
moving linearly towards each other in a way such that their mean remains constant.)
Therefore, to obtain an asymptotic upper bound on $B^{t/n} \Phi_2^{(C_2)}(p^n)$ as $t\downarrow 0$,
we only need to consider the event when exactly one job arrival occurs in the interval $[0,t/n]$,
when the process in some state $q^n$, with the expected change of $\Phi_2^{(C_2)}$ being equal to
$\Delta^{(n)}_{C_2}(q^n)$. But, $q^n \to p^n$ as $t\downarrow 0$, uniformly on all initial points $p^n$ 
within the compact set $\Phi_2(p^n) \le C_2$, and recall that $\Delta^{(n)}_{C_2}(p^n)$ is continuous.
The result follows. 
$\Box$

\subsubsection{Proof of \eqn{eq-generator-bound}.}
\label{sec-generator-bound}

Denote by $\Delta^{(n)}(p^n)$ the expected change of $\Phi_2(p^n)$ when a job arrival occurs at state $p^n$.
Consider a given $p^n$ such that $\Phi_2(p^n) \le C_2$. Then,
obviously, $\Delta^{(n)}_{C_2}(p^n) \le \Delta^{(n)}(p^n)$.
Denote by $\zeta_i=\zeta_i^{(\mathring x^n)}$ the (random) displacement of $w_i$ due to a job arrival.
Then, 
$$
\Delta^{(n)}(p^n)= \E \frac{1}{n} \sum_i (w_i+\zeta_i)^2 - \frac{1}{n} \sum_i w_i^2 = 
\E \frac{1}{n} \sum_i [2w_i \zeta_i + \zeta_i^2],
$$
where expectation $\E$ is with respect to the distribution of the job class, selection set and the component sizes.
For $\zeta_i$ we have:
$$
\zeta_i = \kappa_i - \frac{1}{n} \sum_s \kappa_s, 
$$
where $\kappa_i=\kappa_i^{(\mathring x^n)}$ the (random) actual jump size of particle $i$ due to a job arrival.
Note that $\E \zeta_i$ is the quantity of the order $O(1/n)$, since $\sum_s \kappa_s$ is of order $O(1)$
and $\E \kappa_i $ is of order $O(1/n)$ (because $\P\{i\in S\} = O(1/n)$). Therefore, $\bar \zeta_i = n \E \zeta_i = O(1)$,
and we can write:
$$
\E \sum_i w_i \zeta_i =  \frac{1}{n} \sum_i w_i \bar \zeta_i.
$$
Next,
$$
\E \sum_i \zeta_i^2 \le 2 \sum_i  \E \kappa_i^2  +
2n \frac{1}{n^2} \E (\sum_s \kappa_s)^2 \le C'_1,
$$
because $\sum_i  \E \kappa_i^2$ is uniformly upper bounded by the maximum (job among classes) second moment 
of the total size of all components. Assembling these bounds, we obtain
\beql{eq-n-times-drift}
n\Delta^{(n)}(p^n) 
\le 2 \frac{1}{n} \sum_i w_i \bar \zeta_i + C'_1,
\eeql
and recall that $\Delta^{(n)}_{C_2}(p^n) \le \Delta^{(n)}(p^n)$. This completes the proof. 
$\Box$

\subsection{Proof of Lemma~\ref{lem-small-sufficient}.}
\label{sec-proof-small-sufficient}

Consider the class $j$ as in the lemma statement. WLOG, we can assume that this class is such that 
$m_j$ components are zero and the remaining $d_j - m_j$ components distribution is such that those components are i.i.d.,
uniformly distributed in $[0,\epsilon]^{d_j - m_j}$. (We can always split the original class $j$ into two, with one of the new classes
satisfying the condition above.)
 Fix any $C>0$ as in the definition of a small set. Fix arbitrarily small $\delta_2>0$.
Then, uniformly in all initial states with norm at most $C$, with strictly positive probability, the system gets into a state with norm (= max distance between particles) at most $\delta_2$. Now consider the special system state such that all particles are exactly co-located.  Further, consider a random state, obtained by a sufficiently large, fixed number (say $n$) of  arrivals of this class, such that each particle has non-zero component at least once, with the 
included additional condition that all generated non-zero components happen to fall into $[\epsilon/3, 2\epsilon/3]$; the distribution of this random state -- let us denote it $\mathcal D$ -- is non-proper - the total measure is less than $1$.
We claim the following: if $\delta_2$ is small enough, and the initial state norm is at most $\delta_2$, and we generate a new state 
exactly the same way as we generated $\mathcal D$, except the non-zero components are distributed in $[0,\epsilon]$ (i.e., without including the
additional condition), then the distribution of new random state will dominate $\mathcal D$. The claim is true, because for each outcome in the special system, there exists an ``equally probable'' outcome in the actual system, obtained by slightly shifting, for each particle, its first the non-zero component size to compensate for the deviation of the initial state from special initial state. Therefore, the system satisfies Assumption~\ref{assump-small}. 
$\Box$

\subsection{Proof of Theorem~\ref{thm-bounded-first-moment-tws}.}
\label{sec-proof-bounded-first-moment}

Essentially, this proof uses simplified versions of the arguments used 
in the proof of Theorem~\ref{thm-bounded-first-moment}.

All  ML-FPs $x$ in the family are uniformly Lipschitz -- this follows from \eqn{eq-de-fp-new}.
Each function $x_w$ is differentiable a.e. (with left and right derivatives existing everywhere),
 with strictly negative derivative.

First, we need to show that for each such $x$, the mean $\bar x$ is well-defined finite. We have the characterization of $x$ as an OP-FP, i.e. $x$ is a 
stationary distribution of a particle location evolving in field $x$. Let $X(\cdot)$ be this stationary process of the evolving particle. 
Denote by 
$$
\bar \psi_w = -v + \alpha \int_0^\infty \bar J^{(x,w)}(u) du, 
$$ 
the average drift of the particle at point $w$. It is continuous strictly decreasing at all $w$, except maybe a left boundary point (if any),
where it may have a discontinuity -- jump down. The average of $\bar \psi_w$ w.r.t. distribution $x$ is $0$.
This means that $\bar \psi_w < 0$ for some $w$ sufficiently far on the right. 

If there is a left boundary, we can use $w^2 \wedge C$
as a Lyapunov function. Using the fact that the second moment of the particle jump size 
is uniformly upper bounded by some constant $C_2$, we can derive the ``generator upper bound'' 
(compare to Lemma~\ref{lem-generator-bound})
$$
A_{C_2}(w) \doteq
    \begin{cases}
       2 w \bar \psi_w +C, & \text{if } w^2 \le C_2 \\
      0, & \text{otherwise}
    \end{cases}
$$
for some constant $C>0$, in the sense that
$$
\limsup_{t \downarrow 0} \sup_{w}     \left[\frac{P^t - I}{t} (w^2 \wedge C_2) -  A_{C_2}(w)\right] \le 0.
$$
From here we obtain 
\beql{eq-steady-drift-tws}
0 = \lim_{t\downarrow 0} \E \frac{(X(t)^2 \wedge C_2) - (X(0)^2 \wedge C_2)}{t} \le \E A_{C_2} (X(\infty)) = \E \left[[2 X(\infty) \bar \psi_{X(\infty)} +C] \bI\{X(\infty)^2 \le C_2\}\right],
\eeql
which implies that $\bar x = \E X(\infty) = \infty$ is impossible, because then the RHS of the last display would be negative for large $C_2$. Therefore, 
$\bar x < \infty$.

If there is no left boundary, $\bar \psi_w$ is continuous, and there is a point $u$ such that $\bar \psi_u = 0$.
Let us choose this point for centering and consider Lyapunov function $(w-u)^2 \wedge C_2$.
Using analogous generator upper bound, we can show that $\E |X(\infty) - u| = \infty$ is impossible, because it would
imply negative average drift of $[X(t)-u]^2 \wedge C_2$ when $C_2$ is large. Therefore, $\E |X(\infty) - u| < \infty$, and $\bar x$ is well-defined finite.
(A small subtlety here. If $X(t)$ jumps so that $X(t+) > u + \sqrt{C_2}$, then replacing $[X(t+)-u]^2 \wedge C_2 = C_2$ by $[X(t+)-u]^2$ in the
drift estimate makes the estimate larger.)

We conclude that $\bar x$ is well-defined finite is all cases, and then $\mathring x$ is well-defined. From now on in this proof, WLOG, we assume that $\bar x=0$, so that 
$\mathring x =x$.

Now, suppose the theorem statement does not hold. Then, there exists a sequence $x^{(\kappa)}$ of ML-FP, indexed by $\kappa \uparrow \infty$,
along which $\Phi_1(x^{(\kappa)}) \uparrow \infty$, and all other system parameters converge
(to a set of parameters within the specified class). In this case, 
obviously, both the left boundary (if any) and right boundary (if any) must go to, respectively, $-\infty$ or $+\infty$.
We now show that there exists $\epsilon>0$, such that for all sufficiently large $\kappa$, 
\beql{eq-drift-tws}
G^{(\kappa)}(x^{(\kappa)}) \doteq \E X^{(\kappa)}(\infty) \bar \psi_{X^{(\kappa)}(\infty)} \equiv 
\int_{-\infty}^{\infty} d[-x^{(\kappa)}_w] w \bar \psi^{(\kappa)}_w \le -\epsilon \Phi_1(x^{(\kappa)}).
\eeql
(Compare to \eqn{eq-drift}.) To prove \eqn{eq-drift-tws} we use the same argument as in the proof of 
\eqn{eq-drift} in 
Theorem~\ref{thm-bounded-first-moment}, except we do {\em not} need Assumption~\ref{cond-ntriv-add-i}
to claim the analog of property \eqn{eq-drift-diff} in the case (a) (see text following \eqn{eq-drift-diff}). This is
because $x^{(\kappa)}_w$ are uniformly Lipschitz and, therefore, for a sufficiently small $\delta_3>0$, 
the $x^{(\kappa)}$-measure of the interval  $[a^+_\nu-\delta_3, a^+_\nu]$ is uniformly small. 

Let us use \eqn{eq-steady-drift-tws}, which has the form
$$ 
0 \le \E \left[ [2 X^{(\kappa)}(\infty) \bar \psi_{X^{(\kappa)}(\infty)} +C] \bI\{X^{(\kappa)}(\infty)^2 \le C_2\}\right].
$$ 
Since, for each $\kappa$, $| \bar \psi_w|$ is uniformly bounded and $\E X^{(\kappa)}(\infty)$ is well-defined finite,
by letting $C_2 \uparrow \infty$ we obtain
$$ 
0 \le \E [2 X^{(\kappa)}(\infty) \bar \psi_{X^{(\kappa)}(\infty)} +C],
$$ 
and then, using \eqn{eq-drift-tws},
$$
0 \le -2\epsilon \Phi_1(x^{(\kappa)}) + C.
$$
This means that $\lim_{\kappa} \Phi_1(x^{(\kappa)}) = \infty$ is impossible. The contradiction completes the proof.
$\Box$

\subsection{Proof of Theorem~\ref{thm-free-wave-unique-for-speed}.}
\label{sec-free-wave-unique-for-speed}

\begin{lem}
\label{lem-exp-left-tail}
For any $v$ satisfying \eqn{eq-cond-left-tail-exists}, for any $\epsilon \in (0,1)$ sufficiently close to $1$, 
the proper ML-MFP $x^*$ for the $(\epsilon,0)$-reduced system with frame $(-\infty,0]$, exists. 
Moreover, for this system and such values of $\epsilon$, 
the left tail of any proper ML-FP is asymptotically exponential with the exponent $\beta$ 
determined by \eqn{eq-cond-beta}; namely,
\beql{eq-left-tail-asymp}
\lim_{w\downarrow -\infty} \frac{- \log [1-x^*_w]}{w} = \beta.
\eeql
\end{lem}

 {\em Proof.}  
 If we choose $\epsilon$ sufficiently close to $1$, then the ML-MFP for the $(\epsilon,0)$-reduced system with the right boundary
 will dominate the ML-MFP for the $(\epsilon,0)$-reduced artificial system
 (introduced just above Theorem~\ref{thm-free-wave-unique-for-speed} statement),
 in which, however, the arrival rate $\lambda$ is a little reduced.
 (The flow of arrivals in such artificial system is dominated by the flow of arrivals in the original system, thinned so that 
 any job arrival selecting more than one server in the reduced system is discarded.) But an ML-FP is invariant w.r.t. 
 rescaling of $\lambda$ and $v$ by the same factor. We conclude that the ML-MFP $x^*$ for the $(\epsilon,0)$-reduced system
 dominates the ML-MFP for the $(\epsilon,0)$-reduced artificial system with a slightly larger speed $v' > v$. Therefore, it exists.
 
 Consider now any proper (on the left) ML-FP $x$ for the $(\epsilon,0)$-reduced system. If $\epsilon$ is sufficiently close to $1$,
 the same argument as just above applies to this ML-FP as well. Therefore, $x$ also dominates 
the ML-MFP $\chi^{(\epsilon)}$ for the artificial system with a slightly larger speed $v' > v$. Indeed, 
here we can use the fact that an ML-FP is also a OP-FP, i.e. it is a fixed point of an operator mapping an ``environment'' into
 a stationary distribution of a particle moving in this environment. In this case, for ML-FP $x$, the particle movement is such that, in addition to the jumps it makes in the lower bounding ML-MFP $\chi^{(\epsilon)}$ (for the artificial  system with larger speed or, equivalently, 
 smaller arrival rates), it makes additional jumps. Both these processes have stationary distributions. So, by simple coupling,
 the ML-MFP $\chi^{(\epsilon)}$ for the artificial system is dominated by the ML-FP $x$.

 Since $v' > v$ can be arbitrarily close to $v$, we also see that
 \beql{eq-left-tail-asymp-lower}
\liminf_{w\downarrow -\infty} \frac{- \log [1-x_w]}{w} \ge \beta.
\eeql
Now, let us consider a sequence of $w^* \downarrow -\infty$, and for each $w^*$ we consider the ML-FP $x$ in interval $(-\infty,w^*]$,
which we shift to the right so that the right boundary is at $0$, and normalize so that the scaled number of servers is $1$
instead of $1-x_{w^*}$. Namely, for each $w^*$, consider 
$$
x^{(w^*)}_w = [x_{w+w^*} - x_{w^*}]/[1-x_{w^*}], ~w \le 0.
$$
If we view these functions as distributions, they uniformly dominate the exponential distribution $\chi^{(0)}(\cdot)$,
which is the ML-MFP for the artificial system with speed $v'$ slightly greater than $v$. Therefore, for any subsequence of $w^*$
there is a further subsequence, along which the distributions converge to a distribution $\chi^*$ dominating
$\chi^{(0)}$. Recall that speed $v'$ of $\chi^{(0)}$ is larger than, but can be chosen arbitrarily close to $v$.
As $w^* \downarrow -\infty$, the arrival  rates (per server) $\{\sigma_j\}$ in the reduced systems converge to those in the artificial system.
Using this it is easy to see that any sub-sequential limit of $\chi^*$ must be exactly equal to $\chi^{(0)}$
for the artificial system with the speed exactly $v$. This implies
 \beql{eq-left-tail-asymp-upper}
\liminf_{w\downarrow -\infty} \frac{- \log [1-x_w]}{w} \le \beta,
\eeql
which together with \eqn{eq-left-tail-asymp-lower} gives \eqn{eq-left-tail-asymp}. $\Box$

\begin{lem}
\label{lem-exp-left-tail-unique}
For any $v$ satisfying \eqn{eq-cond-left-tail-exists}, and any $\epsilon \in (0,1)$ sufficiently close to $1$,
the ML-MFP $x^*$ in Lemma~\ref{lem-exp-left-tail} is the unique ML-FP for the $(\epsilon,0)$-reduced system with the right boundary at $0$.
\end{lem}

 {\em Proof.}  Let $\epsilon<1$ be sufficiently close to $1$, so that the ML-MFP $x^{*,\epsilon}$ exist
 for the $(\epsilon,0)$-reduced system with the right boundary, for speed $v$ and any speed $\hat v > v$ sufficiently close to $v$.
 It is not hard to see that the dependence of the ML-MFP on speed is continuous.
 
 The rest of the proof is by contradiction. Suppose that in addition to the ML-MFP $x^{*,\epsilon}$ for speed $v$, there exists a different ML-FP $x$ for the 
 same speed. ML-FP $x$ is dominated by, but not equal to, $x^{*,\epsilon}$; also, 
 by Lemma~\ref{lem-exp-left-tail}, $x$ has an asymptotically exponential left tail
 with the same exponent $\beta$ as $x^{*,\epsilon}$ does. Let us choose speed $\hat v > v$ sufficiently close to $v$, so 
 that ML-MFP $\hat x^{*,\epsilon}$ for speed $\hat v$ is such that
 in some interval $\hat x^{*,\epsilon}_w > x_w$. 
 The left tail of $\hat  x^{*,\epsilon}$ is asymptotically exponential with exponent $\hat \beta < \beta$;
 this, in particular means that for all $w$ sufficiently far on the left
 (i.e. negative, with sufficiently large $|w|$) $\hat  x^{*,\epsilon}_w < x_w$. 
 Therefore, as $w$ increases from $-\infty$ to $0$, the difference $\hat  x^{*,\epsilon}_w - x_w$ is first negative,
 then there is a point when it is strictly positive, and at the right boundary $0$ it is $0$.
 This situation is impossible by Lemma~\ref{lem-hor-dist-no-local-max-new}. 
 $\Box$

Lemma~\ref{lem-exp-left-tail-unique} says that, for any $v$ satisfying \eqn{eq-cond-left-tail-exists},
for any $\epsilon \in (0,1)$ sufficiently close to $1$,
ML-FP for the $(\epsilon,0)$-reduced system with the right boundary at $0$ is unique. But, by Lemma~\ref{lem-de-sol-unique}, any DE-FP $x$ 
is uniquely determined by its left tail, taken as initial condition. Therefore, as a corollary, we obtain the following

 \begin{lem}
\label{lem-ml-fp-unique}
For any $v$ satisfying \eqn{eq-cond-left-tail-exists}, and any $\epsilon \in [0,1)$, proper
ML-MFP $x^{*,\epsilon}$ (if exists) for the $(\epsilon,0)$-reduced system with frame $(-\infty,0]$, is the unique proper ML-FP  for this system.
\end{lem}

{\em Conclusion of the proof of Theorem~\ref{thm-free-wave-unique-for-speed}}.
We see that, for any $v$ satisfying \eqn{eq-cond-left-tail-exists}, there is a unique function $x^*$ which determines
all ML-FP with proper left tail. Namely, choose $\epsilon'$ sufficiently close to $1$, and consider the corresponding
unique DE-FP $x^{*,\epsilon'}$ in the interval $(-\infty,0]$, such that $x^{*,\epsilon'}_0=\epsilon'$.
(This DE-FP $x^{*,\epsilon'}$ is the rescaled version of the unique ML-FP = ML-MFP for the 
$(\epsilon',0)$-reduced system with frame $(-\infty,0]$.) 
Continue this DE-FP solution to the positive half-axis. 
The solution will either (a) hit the $w$-axis or (b) will converge to some $\epsilon_* \ge 0$ as $w\uparrow \infty$.
In case (a): for any $0 \le \epsilon <1$ there is a unique proper ML-FP for the $(\epsilon,0)$-reduced right-regulated system.
In case (b): for any $\epsilon_* < \epsilon <1$ there is a unique proper ML-FP  for the $(\epsilon,0)$-reduced right-regulated system;
for $\epsilon = \epsilon_*$ there is a unique proper ML-FP the $(\epsilon,0)$-reduced free system;
for $\epsilon < \epsilon_*$ there is a unique ML-FP with proper left tail, but not proper on the right, for the $(\epsilon,0)$-reduced free system.
$\Box$ 

\subsection{Sketch of the proof of Theorem~\ref{thm-ssai-right}}
\label{sec-right-regulated-ssai-proof}

First, the result analogous to Theorem~\ref{thm-bounded-first-moment} can be obtained for the right-regulated system with frame $(-\infty,0]$
as well, with analogous proof. 

Then, it is easy to see that the property of the right-regulated system, which is an analog of Condition~\ref{cond-unique} for the
left-regulated system, holds ``automatically.'' This uniqueness property is as follows: {\em for any $v<v_{min}$, for any sufficiently small $\delta \ge 0$,
there exists the unique proper ML-FP ( = ML-MFP) for the $(0,\delta)$-reduced right-regulated system.} 
To prove this, we can take the proper right-regulated ML-FP for the speed $v' \in (v,v_{min})$,
slightly greater that $v$, and scale down its speed to $v$ and scale down input rates $\{\sigma_j\}$ by the same factor $v/v'$; this will give us
a lower bound on the ML-MFP for the $(0,\delta)$-reduced system for sufficiently small $\delta>0$; the latter ML-MFP therefore exists.

Now, an analog of Proposition~\ref{prop-ssai-key} can be obtained for the right-regulated system, with analogous proof.
But this analog of Proposition~\ref{prop-ssai-key} for the right-regulated system does {\em not} require 
an analog of Condition~\ref{cond-unique}, because it (the above uniqueness property) holds automatically.

If we have an analog of Proposition~\ref{prop-ssai-key}, the proof of Theorem~\ref{thm-ssai-right} is 
analogous to that of Theorem~\ref{thm-free-wave-to-ssai}. 
We note the the notion of the $n$-th prelimit system load, which is used in the proof (and defined naturally for the left-regulated system
as the average fraction of particles at the boundary), for the right-regulated system with speed $v$ is defined as follows. 
Denote by $\gamma^n$ the steady-state expected total size of all components of an arriving job. 
Then the ``load'' $\rho^n$ is defined as $\rho^n = v / (\lambda^n \gamma^n)$.
Using coupling, it is easy to observe that $\lambda^n \gamma^n - v$, is non-increasing in $v$, where $\gamma^n$ also depends on $v$.
This, in particular, implies that load $\rho^n$ is non-decreasing in $v$; the continuity in $v$ is also easy to see.
$\Box$

\section{Discussion}
\label{sec-discussion}

This paper leaves several interesting problems open. The results of the paper do not rule out the case when $v_{min} < v_{max}$. However, we do not have an example when this actually occur. Moreover, the results do not rule out the case when $v_{min} < v_{max}$ and for {\em each} speed $v \in (v_{min}, v_{max}]$, in addition to the unique ML-FP for the free system, we have an {\em entire infinite family} of left-regulated ML-FP. 
Again, we do not 
have an example of this situation actually occurring. While such situation may feel counterintuitive, it does occur {\em in different particle system models},
such as the model in \cite{BaSt2023}. (Note that, as we prove, for the right-regulated system the situation is always simpler -- an ML-FP, if any, is unique for each speed $v$.)

In this paper we do identify one case when ``all possible nice properties hold,'' namely the case when Condition~\ref{cond-unique} holds for any speed $v$.
(In turn, Condition~\ref{cond-unique} necessarily holds
when all job classes are $D$-monotone. And, in turn, if all classes have i.i.d. component sizes with IHR, this is sufficient for $D$-monotonicity.) 
Our results show that in this ``nice'' case: (a) $v_{min} = v_{max} = v_*$ and, therefore, the free ML-FP is unique and has speed $v_*$;
(b) SSAI for the left-regulated system holds for the full range of loads; (c) the average advance velocity $v_n$ of the free system 
(with speed parameter $0$) converges to $v_*$. 
Our results do not rule out that Condition~\ref{cond-unique} for any speed $v$ always holds for our general model, 
and then properties (a)-(c) always hold. Or, it is possible that Condition~\ref{cond-unique} does not always hold,
but the nice properties (a)-(c) still hold always. Addressing these questions may be a subject of future work.

Finally, we note that, as far as mean-field analysis of systems with multi-component jobs in concerned, our analysis is quite generic, and may be applicable to other models as well, not necessarily to c.o.c. mechanism. As an example, it applies to systems with water-filling mechanism,
studied in \cite{SnSt2020}. This mechanism is work-conserving and \cite{SnSt2020}, in particular, proves SSAI for the full range of loads.
But, work-conservation is a special case of $D$-monotonicity. As a result, (a simplified version of) the analysis in this paper can be used 
to show that a particle system with water-filling mechanism (under very mild appropriate assumptions) possesses the nice properties
(a)-(c) above.

%\iffalse
%%%%%%%%%\bibliographystyle{acmtrans-ims}
%%%%%%%%%\bibliographystyle{apt}
%\bibliographystyle{abbrv}

%\bibliography{biblio-stolyar}
%%%%%%%%%%%\bibliography{bibliography}
%\fi

\iffalse

\fi

\appendix

\section{Free system positive recurrence}
\label{app1}

\begin{thm}[Generalization of the positive recurrence result in theorem 11 in \cite{St2021-coc}]
\label{th-free-stabil-gen}
Consider a free system. 
Suppose Assumptions~\ref{cond-ntriv-add-ii} and \ref{assump-small} hold. (Assumption~ \ref{cond-ntriv-add-i} is {\em not} needed.)
Then, for all sufficiently large $n$,  $\mathring x^n(\cdot)$ 
is positive recurrent.
\end{thm}

The proof will, in particular, use reduced systems. For a fixed $n$, the reduced system $\{\ell+1, \ldots, \ell+s\}$, consisting of $s$ particles, is such that the out of the ``remaining'' $n-s$ particles, 
$n-\ell-s$ particles are located at $+\infty$ and $\ell$ particles are at $-\infty$. (So, this is the $((n-\ell-s)/n,\ell/n)$-reduced system, as defined earlier in this paper.)

{\em Proof of Theorem~\ref{th-free-stabil-gen}.} Like in the proof of \cite[theorem 11]{St2021-coc}, 
we use the fluid limit technique. 
Let $n$ (and $\{\sigma_j^n\}$) be fixed. Consider a sequence of processes $p^{n,(r)}(\cdot)$, indexed by $r\uparrow \infty$, with initial states such that 
$(1/r) p^{n,(r)}(0) \to q(0)$. 
A process $q(\cdot)$ is a fluid limit of $p^n(\cdot)$, if it is a distributional (weak) limit of the sequence of processes $p^{n,(r)}(\cdot)$.
When we say that a fluid limit (or its trajectories) satisfy a certain property, we mean that they are satisfied w.p.1.
The fluid limit is {\em stable} if for some fixed $T>0$, any fluid limit trajectory starting from unit-norm initial state, $|q(0)| = 1$,
is such that $|q(t)|=0$ for all $t \ge T$. Given Assumption~\ref{assump-small} holds, stability of the fluid limit is sufficient for stability (positive recurrence)
of $p^{n}(\cdot)$.

Recall that $p^n(\cdot)$ is the centered process, and $q(\cdot)$ is its fluid limit. 
In this proof it will be more convenient to consider the fluid limit $\gamma(t) = (\gamma_1(t), \ldots, \gamma_n(t))$ of non-centered (unanchored) free system.
It is such that $\gamma(0)=q(0)$, and all $\gamma_i(t)$ are uniformly Lipschitz non-decreasing.
Stability of fluid limit $q(\cdot)$ is equivalent to the property that, for some $T>0$, any fluid limit $\gamma(\cdot)$, with $|\gamma(0)| =1$, is such that
\beql{eq-fluid-nocenter-stabil}
\gamma_i(t) = \hat \gamma(t), ~\forall t\ge T, ~\forall i \in I,
\eeql
i.e. all particle locations collapse (in fluid limit) after a finite time. We proceed to prove \eqn{eq-fluid-nocenter-stabil}. 

The derivative $\gamma'_i(t)$, which exists a.e. in time, will be referred to as particle $i$ ``velocity'' (in fluid limit).
A maximal subset of co-located (in fluid limit) particles we will call a {\em batch}. Any fluid limit trajectory is such that,
a.e. in time, all particles within a given batch have equal velocities.
First, observe the following. Suppose a state $\gamma(t)$ at time $t$ is such that none of the particles are co-located, that is
$$
\gamma_1(t) < \gamma_2(t) \ldots < \gamma_{n-1}(t) < \gamma_n(t),
$$ 
where we label the particles in the increasing order of their locations. Then, in some positive-length time interval starting time $t$, each particle $i$ will move forward 
at a well-defined constant velocity $u_i \ge 0$, and $u_i$ is non-increasing in $i$ (by monotonicity). If at a given time particles $\ell+1, \ldots, \ell+s$ form a batch (with $\ell$ other particles located strictly behind them and $n-(\ell+s)$ particles strictly ahead of them), and if they have the same common velocity, 
this velocity must be within the interval $[u_{\ell+s},u_{\ell+1}]$ (by monotonicity). From here it follows that, starting any initial state, any fluid limit trajectory is such that if any two particles become co-located at some time, then they will stay co-located from that time on; in other words, the batches cannot decrease in size -- they can only grow by merging with each other. 

Now, suppose we would have a special case that all ``individual'' velocities $u_i$ are distinct: $u_1 > u_2 > \ldots > u_{n-1} > u_n$. 
In this case, given the above observations, any two distinct batches, at any time, have different velocities (the one behind has strictly larger velocity), with the difference bounded away from $0$. This would immediately imply \eqn{eq-fluid-nocenter-stabil}, because all batches would merge into one within finite time. 
In general, however, 
it is easy to see that it is possible that several $u_i$, which are close -- within $\bar d$ -- to the boundary indices $1$ or $n$, may in fact be equal; while the rest of $u_i$ are distinct. 
(As an example, think of a system with a single class $j$ with $(d_j,k_j)=(5,3)$. Then, $u_1=u_2=u_3$ and $u_{n-1}=u_n$.)
Specifically, the general situation is such that
\beql{eq-indiv-vel-order}
u_1  \ge \ldots \ge u_{\bar d} > u_{\bar d +1} > \ldots > u_{n-\bar d} \ge \ldots \ge u_n. 
\eeql
We see that, within a finite time, a batch of particles $\ell +1, \ldots, \ell+s$ with $\ell \le \bar d$ and $\ell+s \ge n-\bar d$
 will form (where particles are indexed in the order of increasing location); in other words, a batch will form, which contains all particles $\bar d +1, \ldots, n-\bar d$. 
Then, the proof of the theorem is completed by the following
\begin{lem}
\label{lem-fluid-key}
A fluid-limit trajectory is such that, a.e. in time $t$, if there is
a batch of particles $\ell +1, \ldots, \ell+s$ of size $s \ge 2\bar d -1$ at time $t$, with the (common) velocity $u_{\ell+1,\ell+s}(t)$, 
then for some fixed $\epsilon'>0$:\\
(i) $u_{\ell+1,\ell+s}(t) \ge u_{\ell+s+1} +\epsilon'$ (if $\ell+s < n$) and (ii) $u_{\ell+1,\ell+s}(t) \le u_{\ell}-\epsilon'$ (if $\ell >0$). \\
(It can be shown that the batch velocity $u_{\ell+1,\ell+s}(t)$ is, in fact, constant $u_{\ell+1,\ell+s}$, which is well-defined in terms of 
the stationary distribution of the reduced system $\{\ell +1, \ldots, \ell+s\}$. We do not need this stronger property for the proof of Theorem~\ref{th-free-stabil-gen}.)
\end{lem}

Indeed,  we already established that a batch containing $\{\bar d+1, \ldots, n-\bar d\}$ will form within a finite time. If $n$ is large, then 
then the size of this batch is large, so that Lemma~\ref{lem-fluid-key} applies.
Then, repeated use of Lemma~\ref{lem-fluid-key} shows that all particles will merge into a single batch within a finite time. 
The proof of Lemma~\ref{lem-fluid-key} is given below.
$\Box$

{\em Proof of Lemma~\ref{lem-fluid-key}.} Let $j$ be a fixed job class, for which Assumption~\ref{cond-ntriv-add-ii} holds. 

(i) $u_{\ell+1,\ell+s}(t) \ge u_{\ell+s+1} +\epsilon_1$ (if $\ell+s < n$). Consider a prelimit system (not the fluid limit), in which we consider
the joint evolution of the reduced system $\{\ell+1, \ldots, \ell+s\}$ and 
the reduced system consisting of a single particle $\ell+s+1$ located infinitely far ahead. It will suffice to show the following: 
upon a class-$j$ job arrival, regardless of the current state of the reduced system, the 
distribution of the random advance of a particle in the reduced system, picked uniformly at random, and the distribution of the advance of
particle $\ell+s+1$ can be coupled so that with a positive probability at least $\delta_1>0$, the former advance is greater than the latter by at least  
by some $\epsilon_1>0$. (In the fluid limit, using strong law of large numbers, this translates into an instantaneous velocity of the batch being 
greater than the constant velocity of particle $\ell+s+1$ at least by some $\epsilon'>0$. 
We also use the fact that a finite distance between batches on a fluid limit translates into an increasing to infinity distance in a pre-limit system.)
We proceed with proving this fact.

Observe that, with probability at least $1/s$, the picked particle $i$ in the reduced system $\{\ell+1, \ldots, \ell+s\}$
is a left-most particle in it (in the pre-limit process). Then, it suffices to couple and compare the advance of $i$ under the condition 
that particle $i$ is selected by the job and particle $\ell+s+1$ is not, to the advance of $\ell+s+1$ under the condition 
that particle $\ell+s+1$ is selected and particle $i$ is not, and show that the former advance is larger with positive probability.
(This is because, using monotonicity, we can always couple {\em any} two selection events, which are ``symmetric'' with respect to particles $i$ and $\ell+s+1$, so that 
the advance of $i$ dominates that of $\ell+s+1$.)
The latter property is, in fact, true. Indeed, consider the natural coupling of those two events, namely
such that particles $i$ and $\ell+s+1$ receive exactly same job component size, and the remaining selected servers and their component sizes are common.
Note that, under this coupling, particle $\ell+s+1$ can never receive a larger displacement than $i$. Further, for some fixed $\epsilon > 0$, the following event occurs with probability at least some $\delta>0$: the other $d_j-1$ selected particles (besides $i$ or $\ell+s+1$) are the particles $i+1, \ldots, i+d_j-1$ (in the order of increasing location), with their component sizes being non-decreasing in location and such that $\xi_{i+k_j} \ge \epsilon$, and   $\xi_{i}\ge \epsilon$ [or, respectively, 
$\xi_{\ell+s+1} \ge \epsilon$]. Under such event, particle $i$ advances at least by $\epsilon>0$, while particle $\ell+s+1$ does not advance at all. This concludes the proof of (i).

(ii) $u_{\ell+1,\ell+s}(t) \le u_{\ell}-\epsilon'$ (if $\ell >0$). The argument is similar to the one we used to prove (i), but here we need to consider 
the evolution of the prelimit system over  {\em two consecutive} class-$j$ job arrivals. Consider a prelimit system (not the fluid limit), in which we consider
the joint evolution of the reduced system $\{\ell+1, \ldots, \ell+s\}$ and 
the reduced system consisting of single particle $\ell$ located infinitely far behind.
It will suffice to show the following: 
upon two consecutive class-$j$ job arrivals, regardless of the current state of the reduced system, the 
distribution of the random advance of a particle in the reduced system, picked uniformly at random, and the distribution of the advance of
particle $\ell$ can be coupled so that with a positive probability at least $\delta_1>0$, the former advance is less than the latter by at least  
by some $\epsilon_1>0$.

Observe that, with probability at least $1/s$, the picked particle $i$ in the reduced system $\{\ell+1, \ldots, \ell+s\}$
is a right-most particle in it (in the pre-limit process). Then, it suffices to show that, under appropriate coupling, with positive probability,
 the advance of $i$ under the condition 
that particle $i$ is selected by both jobs and particle $\ell$ is selected by neither, is strictly smaller than the advance of $\ell$ under the condition 
that particle $\ell$ is selected by both jobs and particle $i$ is selected by neither. Consider the natural coupling of these two events,
with particles $i$ and $\ell$ receiving exactly same component sizes, and the remaining selected servers and their component sizes being common.
There exist some fixed $\epsilon >0$, such that the following event occurs with probability at least some $\delta>0$:
 the first job, besides particle $i$ [resp., particle $\ell$], selects an arbitrary fixed set $\hat I_1$ of $d_j-1$ particles in the reduced system 
 $\{\ell+1, \ldots, \ell+s\}$; the $k_j$-th smallest component of the first job received by particles in $\hat I_1$ has size at least $\epsilon$,
 while the component size received by particle $i$ [resp., particle $\ell$] is the largest of all components of the job (so it it is also at least $\epsilon$);
 the second job has the same properties as the first one, except, besides particle $i$ [resp., particle $\ell$], it selects an arbitrary fixed set $\hat I_2$ of $d_j-1$ particles within  $\{\ell+1, \ldots, \ell+s\} \setminus \hat I_1$. Under this event, the following must hold: either\\
 {\em upon the first job arrival 
 particle $\ell$ advances by at least $\epsilon$ while particle $i$ advances by less than $\epsilon/2$} \\
 or \\
 {\em upon the first job arrival 
 particle $i$ advances by at least $\epsilon/2$ and then after the second job arrival particle $i$ advances by at least $\epsilon/2$ less than 
 particle $\ell$}. \\
 In either case, after the two job arrivals, the advance of particle $\ell$ is at least by $\epsilon/2$ greater than that of particle $i$.
$\Box$

\end{document}